\documentclass[10pt]{amsart}
\usepackage{amsthm,amsmath,amssymb,amscd}
\usepackage{amssymb}
\usepackage{graphics}
\usepackage{color}

\newtheorem{Theorem}{Theorem}[section]      
\newtheorem{Proposition}[Theorem]{Proposition}    
\newtheorem{Lemma}[Theorem]{Lemma}            
\newtheorem{Definition}[Theorem]{Definition}

\title[Extremal domains of big volume in a compact manifold]{Extremal domains of big volume for the first eigenvalue of the Laplace-Beltrami operator in a compact manifold}

\author[Pieralberto Sicbaldi]{Pieralberto Sicbaldi}
\address{pieralberto.sicbaldi@univ-paris12.fr  \\ Universit\'e Paris-Est, France}

\begin{document}

\maketitle
{\bf Abstract.} We prove the existence of extremal domains for the first eigenvalue of the Laplace-Beltrami operator in some compact Riemannian manifolds of dimension $n \geq 2$, with volume close to the volume
of the manifold. If the first (positive) eigenfunction $\phi_0$ of the Laplace-Beltrami operator over the manifold is a nonconstant function, these domains are close to
the complement of geodesic balls of small radius whose center is close to the point where
$\phi_0$ attains its maximum. If $\phi_0$ is a constant function and $n \geq 4$,
these domains are close to the complement of geodesic balls of small radius whose center is close to a nondegenerate critical point of the scalar curvature function.

\section{Statement of the result}

Article [\ref{Pac-Sic}] is a first study on the possibility to construct extremal domains for the first eigenvalue of the Laplace-Beltrami operator in a Riemannian manifold. In that work are introduced all basic definitions and properties of extremal domains for the first eigenvalue of the Laplace-Beltrami operator, and the examples of domains obtained have small volume. In this article we will give an existence result for extremal domains of great volume in a compact riemannian manifold. We will be interested in domains obtained by taking the complement of small domains contained in the interior of the manifold: $M \backslash \Omega$, where $\Omega$ has a small volume and $\overline{\Omega} \subseteq \mathring M$. For the sake of completeness, we start recalling the basic facts of extremal domains (for other details see [\ref{Pac-Sic}]). 

\medskip

Assume that we are given $(M,g)$ a compact $n$-dimensional Riemannian manifold, $n \geq 2$, with or without boundary $\partial M$. In the case $\partial M \neq \emptyset$, then $\partial M$ is supposing to be an $n-1$-dimensional Riemannian manifold. Let $\overline \Omega_0$ be a domain in the interior of $M$ and let us consider the domain $M \backslash \Omega_0$. 
\begin{Definition}
We say that $\{M \backslash \Omega_t \}_{t \in (-t_0, t_0)}$, $\overline{\Omega}_t \subseteq \mathring M$, is a deformation of $M \backslash \Omega_0$ if there exists a vector field $\Xi$ (such that $\Xi(\partial M) \subseteq \partial M$) for which
$M \backslash \Omega_t = \xi (t,M \backslash \Omega_0)$ where  $\xi(t, \cdot)$ is the flow associated to $\Xi$, namely 
\[
\frac{d\xi}{dt} (t,p)= \Xi (\xi(t,p)) \qquad \mbox{and} \qquad \xi(0, p) =p \,.
\] 
The deformation is said to be volume preserving if the volume of $M \backslash \Omega_t$ does not depend on $t$.
\end{Definition}

Let us denote by $\lambda_t$ the first eigenvalue of $-\Delta_{g}$ on $M \backslash \Omega_{t}$ with $0$ Dirichlet boundary condition on $\partial \Omega_t$. In the case where $\partial M \neq \emptyset$, then we ask also one of the following boundary condition :
\begin{enumerate}
	\item $0$ Dirichlet boundary condition on $\partial M$, or
	\item $0$ Neumann boundary condition on $\partial M$.
\end{enumerate}
We will suppose the regularity of $\partial M$.
Observe that both $t \mapsto \lambda_t$ and the associated eigenfunction $t \mapsto u_t$  (normalized to have $L^2 (M \backslash \Omega_t)$ norm equal to $1$) are continuously differentiable. 

\medskip

\begin{Definition}
A domain $M \backslash \Omega_{0}$ is an \textit{extremal domain} for the first eigenvalue of $-\Delta_{g}$ if for any volume preserving deformation $\{{M \backslash \Omega}_t\}_t$ of $M \backslash {\Omega}_{0}$, we have 
\[
\frac{d \lambda_t}{dt} |_{t =0} = 0 \, .
\]  
\label{def:1}
\end{Definition}

According to the condition taken at the boundary (if the boundary is not empty) we will talk about extremal domains under the 0 Dirichlet boundary condition at $\partial M$ or extremal domains under the 0 Neumann boundary condition at $\partial M$. Let $\phi_0$ be the first eigenfunction of the Laplace-Beltrami operator over the manifold M
\begin{equation}
\begin{array}{rcccl}
	\Delta_{g} \, \phi_0 + \lambda_0 \, \phi_0 & = & 0 & \textnormal{in} & M 
\end{array}
\end{equation}
with 0 Dirichlet or 0 Neumann boundary condition (if $\partial M \neq \emptyset$), normalized to have $L^2$-norm equal to 1. Here $\lambda_0$ is the first eigenvalue of $-\Delta_g$ on $M$. If the volume of $\Omega$ is very small, it is natural to expect that the first eigenfunction of the Laplace-Beltrami operator over $M \setminus \Omega$ will be close to $\phi_0$. We remark that we have to distinguish two cases of behaviour of $\phi_0$ (and then also of the first eigenfunction over $M \setminus \Omega$), according with the condition at the boundary :

\medskip

\begin{itemize}
	\item CASE 1. If $\partial M \neq \emptyset$ and $\phi_0$ satisfy the 0 Dirichlet condition on $\partial M$ then $\phi_0$ is a positive non constant function, and then attains its maximum in at least a point of the manifold, say at $p_0$. Moreover $\lambda_0 > 0$.
		
\medskip
	\item CASE 2. If $\partial M = \emptyset$, or if $\partial M \neq \emptyset$ and $\phi_0$ satisfy the 0 Neumann condition on $\partial M$, then $\phi_0$ is a constant function 
	\[
	\phi_0 = \frac{1}{\sqrt{\textnormal{Vol}_g(M)}}.
	\]
and $\lambda_0 = 0$.
	
\end{itemize}
When we will consider the first eigenfunction of the Laplace-Beltrami operator over $M \setminus \Omega$, where $\Omega \subset \mathring M$, we will take at $\partial M$ the same boundary condition of $\phi_0$, distinguishing always the two cases.

\medskip

For all $\epsilon >0$ small enough, we denote by $B_{\epsilon}(p) \subset M$ the geodesic ball of center $p \in M$ and radius $\epsilon$. We denote by  $\mathring B_\epsilon \subset \mathbb R^n$ the Euclidean ball of radius $\epsilon$ centered at the origin. 

\medskip

Now we can state the main result of our paper~:
\begin{Theorem}\label{maintheorem}
In the CASE 1 assume that $p_0$ is a nondegenerate critical point of the first eigenfunction $\phi_0$ of the Laplace-Beltrami operator over $M$, and in the CASE 2 assume that $p_{0}$ is a nondegenerate critical point of $\textnormal{Scal}$, the scalar curvature function of $(M,g)$. In the CASE 2 we will assume also $n \geq 4$. Then, for all $\epsilon > 0$ small enough, say $\epsilon \in (0, \epsilon_0)$, there exists a smooth domain $\Omega_\epsilon \subset M$ such that~:
\begin{itemize}
\item[(i)] The volume of $\Omega_\epsilon$ is equal to the Euclidean volume of $\mathring B_\epsilon$. 
\item[(ii)] The domain $M \backslash \Omega_\epsilon$ is extremal in the sense of definition~\ref{def:1}.
\end{itemize}
Moreover there exist a constant $c > 0$ and for all $\epsilon \in (0, \epsilon_0)$ there exists $p_\epsilon \in M$ such that the boundary of $\Omega_\epsilon$ is a normal graph over $\partial B_\epsilon (p_\epsilon)$ for some function $w_\epsilon$, with
\[
 \textnormal{ dist} (p_\epsilon, p_0) \leq c \, \epsilon \, .
\]
and
\[
\begin{array}{ll}
\|w_\epsilon \|_{\mathcal C^{2, \alpha}\partial B_\epsilon (p_\epsilon))} \leq c \, \epsilon^2 & \textnormal{in the CASE 1 and}\,\,n \geq 3\\[3mm]
\|w_\epsilon \|_{\mathcal C^{2, \alpha}\partial B_\epsilon (p_\epsilon))} \leq c \, \epsilon^2\, \log \epsilon & \textnormal{in the CASE 1 and}\,\,n = 2\\[3mm]
\|w_\epsilon \|_{\mathcal C^{2, \alpha}\partial B_\epsilon (p_\epsilon))} \leq c \, \epsilon^3 & \textnormal{in the CASE 2 and}\,\,n\geq5\textnormal{}\\[3mm]
\|w_\epsilon \|_{\mathcal C^{2, \alpha}\partial B_\epsilon (p_\epsilon))} \leq c \, \epsilon^3\, \log \epsilon & \textnormal{in the CASE 2 and}\,\,n = 4
\end{array}
\]
\label{th:1}
\end{Theorem}

\medskip

We remark that the theorem do not give any information in the CASE 2 for the dimensions 2 and 3. In fact our methode to prove the main theorem is based, for the CASE 2, on the approximation of some Green function to the first eigenfunction of the Laplace-Beltrami operator outside a perturbed ball. When the dimension of $M$ is at least 4, we are able to compute the first coefficients of the local expansion of that Green function and this allows us to obtain the estimations we need. But for the dimensions 2 and 3, global terms depending on the manifold do not allow us to obtain a local expansion of that Green function. This is the reason for which we didn't obtain information on extremal domains of big volume in the CASE 2 for the dimensions 2 and 3. 

\section{Characterization of the problem}

In order to prove our theorem we need the following result that caracterizes extremal domains of the form $M \backslash \Omega$ in a Riemannian manifold $M$. The following result gives a formula for the first variation of the first eigenvalue  for some mixte problems under variations of the domain. 

\medskip
 
We have the~:
\begin{Proposition}\label{lambda}
The derivative of $t \longmapsto \lambda_t$ at $t =0$ is given by
\[
\frac{d\lambda_t}{dt} |_{t =0}= - \int_{\partial \Omega_{0}}  \left( g(\nabla  u_0 ,  \nu_0) \right)^{2}\  g(\Xi, \nu_0) \, \mbox{\textnormal{dvol}}_g,
\]
where $\mbox{\textnormal{dvol}}_g$ is the volume element on $\partial \Omega_{0}$ for the metric induced by $g$ and $\nu_0$ is the normal vector field about $\partial \Omega_0$.
\end{Proposition}
{\bf Proof :} We denote by $\xi$ the flow associated to $\Xi$.  By definition, we have
\begin{equation}\label{eq:1-1}
u_t (\xi (t , p)) = 0
\end{equation}
for all $p \in \partial \Omega_0$. Moreover, if we take the 0 Dirichlet boundary condition on $\partial M$ then equation (\ref{eq:1-1}) is valid also on $\partial M$. On the other hand, if we take the 0 Neumann condition on $\partial M$ then we have
\begin{equation}\label{aa}
g(\nabla u_t (\xi (t , p)), \nu_{t}) = 0
\end{equation}
for all $p \in \partial M$, where $\nu_{t}$ is the unit normal vector about $\partial M$. 

\medskip

Differentiating (\ref{eq:1-1}) with respect to $t$ and evaluating the result at $t = 0$ we obtain  
\[
\partial_t u_0 = -  g (\nabla  u_0 ,  \Xi  ) \, ,
\]
on $\partial \Omega_0$. Now $u_0 \equiv 0$ on $\partial \Omega_0$, and hence only the normal component of  $\Xi$ plays a r\^ole in this formula. Therefore, we have 
\begin{equation}\label{eq:1-2}
\partial_t u_0 = -  \,  g(\nabla  u_0 ,  \nu_0) \, g(\Xi, \nu_0) \, ,
\end{equation}
on $\partial \Omega_0$. The same reasoning is also valid on $\partial M$ if we take the 0 Dirichlet boundary condition on $\partial M$. In this case, by the fact that $\Xi(\partial M) \subseteq T(\partial M)$ we have
\begin{equation}\label{eq:1-2bis}
\partial_t u_0 = 0 \, 
\end{equation}
on $\partial M$. On the other hand, if we take the 0 Neumann condition on $\partial M$ then taking a system of coordinates $x = (x_1, ..., x_n)$ such that $\nu_{t} = -\partial_{x^0}$ on $\partial M$ and differentiating (\ref{aa}) with respect to $t$ and evaluating the result at $t = 0$ we obtain  
\begin{equation}\label{bb}
0 = - \partial_{x^0} \partial_t u_0 - g (\nabla \partial_{x^0} u_{0}, \Xi) = - \partial_{x^0} \partial_t u_0 = g (\nabla \partial_t u_0, \nu_{0})  \, 
\end{equation}
on $\partial M$, where we used the fact that $\nu_{t}$ don't depend on $t$ on $\partial M$ together with the facts that $\partial_{x^0} u_{0} = 0$ on $\partial M$ and that $g (\Xi, \nu_{0} )= 0$ in $\partial M$ because $\Xi(\partial M) \subseteq T(\partial M)$. 

\medskip

We differentiate now with respect to $t$ the identity 
\begin{equation}\label{eq:1-3}
\Delta_g \, u_t + \lambda_t \, u_t  =0.
\end{equation}
and again evaluate the result at $t = 0$. We obtain
\begin{equation}\label{eq:1-4}
\Delta_{g} \partial_t  u_0 + \lambda_0 \,  \partial_t  u_0  = - \partial_t  \lambda_0  \, u_0 \, ,
\end{equation}
in $\Omega_{0}$.  Now we multiply (\ref{eq:1-4}) by $u_0$ and  (\ref{eq:1-3}), evaluated the result at $t =0$, by $\partial_t  u_0$, subtract the results  and integrate it over $\Omega_{0}$ to get~:
\begin{eqnarray*}
\partial_t  \lambda_0 \, \int_{\Omega_{0}} u^2_0 \, \textnormal{dvol}_g & = & \int_{M \backslash \Omega_{0}} \left(\partial_t  u_0 \, \Delta_{g} u_0 - u_0  \, \Delta_{g} \partial_t  u_0 \right) \textnormal{dvol}_g \\[3mm] 
& = & \int_{\partial M \cup \partial \Omega_{0}} \left( \partial_t  u_0 \, g(\nabla  u_0 , \nu_0) - u_0 \,   g(\nabla  \partial_t  u_0 ,  \nu_0 ) \right) \textnormal{dvol}_g \\[3mm]
& = & \int_{\partial \Omega_{0}} \left( \partial_t  u_0 \, g(\nabla  u_0 , \nu_0) - u_0 \,   g(\nabla  \partial_t  u_0 ,  \nu_0 ) \right) \textnormal{dvol}_g \\
&  & + \int_{\partial M} \left( \partial_t  u_0 \, g(\nabla  u_0 , \nu_0) - u_0 \,   g(\nabla  \partial_t  u_0 ,  \nu_0 ) \right) \textnormal{dvol}_g \\[3mm]
& = & - \int_{\partial \Omega_{0}} \left(g(\nabla  u_0 , \nu_0)  \right)^{2}\  g(\Xi, \nu_0) \ \textnormal{dvol}_g \, ,
\end{eqnarray*}
where we have used (\ref{eq:1-2}), (\ref{eq:1-2bis}) or (\ref{bb}), the fact that $u_0=0$ on $\partial \Omega_0$, and the the fact that $u_0=0$ or $g(\nabla u_0, \nu_{0}) = 0$ on $\partial M$ to obtain the last equality. The result follows at once from the fact that $u_0$ is normalized to have $L^2(\Omega_0)$ norm equal to $1$. Observe that in the previous argument $\partial M$ can be empty. \hfill  $\Box$ 

\medskip

This result allows us to characterize extremal domains for the first eigenvalue of the Laplace-Beltrami operator under some particular 0 mixte boundary conditions, and state the problem of finding extremal domains into the solvability of an over-determined elliptic problem. The proof of the following proposition is a consequence of the previous result; because it is very similar to the proof of Proposition 2.2 in [\ref{Pac-Sic}] we don't report it here.

\begin{Proposition}
Given a smooth domain $\Omega_0$ contained in the interior of $M$, the domain $M \backslash \Omega_0$ is extremal if and only if there exists a constant $\lambda_0$ and a positive function $u_0$ (if $\partial M \neq \emptyset$ then we take 0 Dirichlet (CASE 1) or 0 Neumann (CASE 2) boundary condition on $\partial M$) such that
\begin{equation}
\left\{
\begin{array}{rclll}
	\displaystyle \Delta_{g} u_0 + \lambda_0 \, u_0 & = & 0 & \textnormal{in} & M \backslash \Omega_0 \\[3mm]
	\displaystyle u_0 & = & 0 & \textnormal{on} & \partial \Omega_0 \\[3mm]
	\displaystyle g (\nabla  u_0 , \nu_0) & = & \textnormal{constant} & \textnormal{on} & \partial\Omega_0  \, ,
\end{array}
\right.
\label{eq:1-5}
\end{equation}
where $\nu_0$ is the normal vector field about $\partial \Omega_0$.
\end{Proposition}

\medskip

Therefore, in order to find extremal domains, it is enough to find a domain $M \backslash \Omega_0$ (regular enough) for which the over-determined problem (\ref{eq:1-5}) has a nontrivial positive solution. In this article we will solve this problem to find solutions whose volumes are close to the volume of our compact manifold.

\section{Rephrasing the problem}

Following the approach of [\ref{Pac-Sic}], we introduce the following notation.  Given a point $p\in M$ we denote by $E_{1} ,\ldots, E_n$ an orthonormal basis of the tangent plane to $M$ at $p$. Geodesic normal coordinates $x : =(x^1, \ldots, x^n) \in \mathbb R^n$ at $p$ are defined by 
\[
X (x) : =\mbox{Exp}_p^g \left( \sum_{j=1}^n x^j \, E_j \right) 
\]
We recall the Taylor expansion of the coefficients $g_{ij}$ of the metric $X^* g$ in these coordinates.
\begin{Proposition}
At the point of coordinate $x$, the following expansion holds~:
\begin{equation}
g_{ij} = \delta_{ij} + \frac{1}{3} \, \sum_{k,\ell} R_{ikj\ell} \, x^{k} \, x^{\ell} + \frac{1}{6} \,\sum_{k, \ell, m}  R_{ikjl,m} \, x^{k} \, x^{\ell} \, x^{m} + {\mathcal O}(|x|^{4}),
\end{equation}
Here $R$ is the curvature tensor  of $g$ and
\begin{eqnarray*}
R_{ikj\ell} & = & g\big( R(E_{i}, E_{k}) \, E_{j} ,E_{\ell}\big)\\[3mm]
R_{ikj\ell ; m} & = & g\big( \nabla_{E_{m}} R (E_{i}, E_{k}) \, E_{j} , E_{\ell} \big) \, ,
\end{eqnarray*}
are evaluated at the point $p$.
\label{pr:1.1}
\end{Proposition}
The proof of this proposition can be found in \cite{willmore} or also in \cite{schoen}.

\medskip

It will be convenient to identify $\mathbb R^n$ with $T_pM$ and $S^{n-1}$ with the unit sphere in $T_pM$. If  $x : = (x^1, \ldots, , x^n) \in \mathbb R^n$, we set
\[
\Theta (x) : = \sum_{i=1}^n x^i \, E_i \in T_pM \, .
\]
Given a continuous function $f : S^{n-1} \longmapsto (0, \infty)$ whose $L^\infty$ norm is small (say less than the cut locus of $p$)  we define 
\[
B_{f}^g (p) : =  \left\{ \mbox{Exp}_p (  \Theta (x)) \qquad : \quad  x \in \mathbb R^{n} \qquad 0 \leq |x|   < f (x/|x|) \right\} \, . 
\]
The superscript $g$ is meant to remind the reader that this definition depends on the metric.

\medskip

Our aim is to show that, for all  $\epsilon >0$ small enough, we can find a point $p\in M$ and a function $v : S^{n-1} \longrightarrow \mathbb R$  such that 
\[
\textnormal{ Vol} \, B_{\epsilon(1+v)}^g(p) =  \epsilon^n \, \textnormal{ Vol} \, \mathring B_1
\]
and the over-determined problem
\begin{equation}\label{a1}
\left\{
\begin{array}{rcccl}
	\Delta_{g} \, \phi + \lambda \, \phi & = & 0 & \textnormal{in} & M \setminus B_{\epsilon(1+v)}^g (p) \\[3mm]
	\phi & = & 0 & \textnormal{on} & \partial B^g_{\epsilon(1+v)}(p) \\[3mm]
	\displaystyle  g( \nabla  \phi ,  \nu) & = & \textnormal{ constant} & \textnormal{on} & \partial B_{\epsilon(1+v)}^g (p) 
\end{array}
\right.
\end{equation}
with 0 Dirichlet (CASE 1) or 0 Neumann (CASE 2) boundary condition on $\partial M$ if $\partial M \neq \emptyset$,
has a nontrivial positive solution, where $\nu$ is the normal vector field about $\partial B_{\epsilon(1+v)}^g(p)$. 

\medskip

Observe that, considering the dilated metric $\bar g : =  \epsilon^{-2} \, g$, the above problem is equivalent to finding a point $p\in M$ and a function $v : S^{n-1} \longrightarrow \mathbb R$  such that 
\[
\textnormal{ Vol} \, B_{1+v}^{\bar g}(p) =  \textnormal{ Vol} \, \mathring B_{1}
\]
and for which the over-determined problem
\begin{equation}\label{a1bis}
\left\{
\begin{array}{rcccl}
	\Delta_{\bar g} \, \bar \phi + \bar \lambda \, \bar \phi & = & 0 & \textnormal{in} & B_{1+v}^{\bar g} (p) \\[3mm]
	\bar \phi & = & 0 & \textnormal{on} & \partial B^{\bar g}_{1+v}(p) \\[3mm]
	\displaystyle  \bar g( \nabla  \bar \phi , \bar \nu ) & = & \textnormal{ constant} & \textnormal{on} & \partial B_{1+v}^{\bar g} (p) 
\end{array}
\right.
\end{equation}
with 0 Dirichlet (CASE 1) or 0 Neumann (CASE 2) boundary condition on $\partial M$ if $\partial M \neq \emptyset$,
has a nontrivial positive solution, where $\bar \nu$ is the normal vector field about $\partial B_{1+v}^{\bar g}(p)$. We can simply consider 
\[
\phi = \bar \phi
\]
(naturally it will not have the norm equal to 1, but depending on $\epsilon$) and 
\[
\lambda = \epsilon^{-2}\, \bar \lambda \, .
\]
In what it follows we will consider sometimes the metric $g$ and sometimes the metric $\bar g$, in order to simplify the computation we will meet.

\section{The first eigenfunction of $-\Delta_g$ outside a small ball}

We remark that the positive solution of the problem 
\begin{equation}{\label{foutball}}
\left\{
\begin{array}{rcccl}
	\Delta_{g} \, \phi_\epsilon + \lambda_\epsilon \, \phi_\epsilon & = & 0 & \textnormal{in} & M \setminus B_{\epsilon}^g (p) \\[3mm]
	\phi_\epsilon & = & 0 & \textnormal{on} & \partial B^g_{\epsilon}(p)
\end{array}
\right.
\end{equation}
with 0 Dirichlet (CASE 1) or 0 Neumann (CASE 2) condition on $\partial M$ if $\partial M \neq \emptyset$, normalized to have $L^2(M \setminus B_{\epsilon}^g (p))$-norm equal to 1, a priori is not known. In this section we will be interested in that solution.

\medskip

Let $p \in M$, let $c$ be a constant, and let $\Gamma_p$ be a Green function over $M$ with respect to the point $p$ defined by 
\begin{equation}\label{greendefinition}
	- \big ( \Delta_{g} \, + \lambda_0 \big)\Gamma_p  = c_n\, \big( \delta_p - \phi_0(p)\, \phi_0\big)\qquad  \textnormal{in}\qquad  M 
\end{equation}
with 0 Dirichlet boundary condition (for the CASE 1) or 0 Neumann boundary condition (for the CASE 2) at $\partial M$ if $\partial M \neq \emptyset$,, and normalization 
\[
\int_{M} \Gamma\, \phi_0 \, \textnormal{ dvol}_{g} = 0,
\]
where  
$\delta_p$ is the Dirac distribution for the manifold $M$ with metric $g$ at the point $p$, i.e.
\[
\int_{M} \delta_p\, f\, \textnormal{ dvol}_{g} = f(p)\textnormal{ for all }f \in C^{\infty}_0(M) \qquad \textnormal{and} \qquad \int_{M} \delta_p\, \textnormal{ dvol}_{g} = 1
\]
 We remark that $\Gamma_p$ exists because
\[
\int_M \left[\delta_p - \phi_0(p)\,\phi_0\right]\, \phi_0\, \textnormal{ dvol}_{g} = 0.
\]
It is easy to check that for each dimension $n$ of the manifold it is possible to chose the constant $c_n$ in order to have the following expansions in local coordinates $x$ of $\Gamma_p$ in a neighborhood of the point $p$ :
\[
\begin{array}{llclcll}
\textnormal{for} & n = 2 & : & \Gamma_p(x) & = & \log|x| + a + \mathring g (b,x) + \mathcal{O}(|x|^{\alpha}) & \forall \alpha < 2\\[3mm]
\textnormal{for} & n = 3 & : & \Gamma_p(x) & = & |x|^{-1} + a' + \mathcal{O}(|x|^{\alpha}) & \forall \alpha < 1\\[3mm]
\textnormal{for} & n = 4 & : & \Gamma_p(x) & = & |x|^{-2} + \mathcal{O}(|x|^{\alpha}) & \forall \alpha < 0\\[3mm]
\textnormal{for} & n \geq 5 & : & \Gamma_p(x) & = & |x|^{2-n} + {\mathcal O}(|x|^{4-n}) 
\end{array}
\]
where $a,a' \in \mathbb{R}$ and $b \in \mathbb{R}^n$, and $\mathring g(\cdot,\cdot)$ is the scalar product in $\mathbb{R}$.

\medskip

It is useful to consider the weighted space $\mathcal{C}^{k,\alpha}_{\nu}(M \setminus \{p\})$, defined as the space of functions in $\mathcal{C}^{k,\alpha}(M \setminus \{p\})$ such that, in the normal geodesic coordinates $x$ around $p$, 
\begin{equation}\label{gammaexpansion}
\begin{array}{l}
\displaystyle \| u\|_{\mathcal{C}^{k,\alpha}_{\nu}(M \setminus \{p\})} := \sup_{\mathring B_{R_0}} |x|^{-\nu}\, |u| + \sup_{\mathring  B_{R_0}} |x|^{1-\nu}\, |\nabla u| + \sup_{\mathring B_{R_0}} |x|^{2-\nu}\, |\nabla^2 u| + \cdots + \\[3mm]
\qquad \displaystyle + \sup_{\mathring B_{R_0}} |x|^{k-\nu}\, |\nabla^k u| + \sup_{0<R\leq R_0}\, \sup_{x,y \in \mathring B_R \setminus \mathring B_{R/2}}\, R^{k+\alpha-\nu}\, \left| \frac{\nabla^k u(x) - \nabla^k u(y)}{|x-y|^\alpha}\right|\, \leq\, \infty.
\end{array}
\end{equation}
where $R_0$ is chosen in order to have the existence of the local coordinates $x \in \mathring B_{R_0}$.

\medskip

Let us consider $\varphi \in \mathcal{C}^{2,\alpha}_m(S^{n-1})$, where $m$ is meant to point out that functions have (euclidean) mean 0 over $S^{n-1}$, and let $H_{\varphi}$ be a bounded harmonic extension of $\varphi$ to $\mathbb{R}^{n} \setminus \mathring B_1$ :
\begin{equation}\label{Hharmonicc}
\left\{
\begin{array}{rclll}
	\Delta_{\mathring g} H_{\varphi} & = & 0 & \textnormal{in} & \mathbb{R}^{n} \setminus \mathring B_1\\[3mm]
	 H_{\varphi} & = & \varphi & \textnormal{on} & \partial \mathring B_1 
\end{array}
\right.
\end{equation}
where $\mathring g$ is the euclidean metric and we identified $\partial \mathring B_1$ with $S^{n-1}$. We have the :

\begin{Lemma}\label{lemmaharmonic}
The following inequality holds :
\[
\|H_{\varphi}(x)\|_{\mathcal{C}^{2,\alpha}_{1-n}(\mathbb{R}^{n} \setminus \mathring B_1)} \leq c\, \| \varphi \|_{\mathcal{C}^{2,\alpha}(S^{n-1})}
\]
for some positive constant $c$. In particular
\[
\lim_{|x| \rightarrow + \infty} H_{\varphi}(x) = 0.
\]
\end{Lemma}

{\bf Proof.} Let us consider
\[
\varphi = \sum_{j = 1}^\infty \varphi_j 
\]
the eigenfunction decomposition of $\varphi$, i.e. 
\begin{equation}\label{decomposition}
\Delta_{S^{n-1}} \varphi_j = -j(n-2+j)\, \varphi_j
\end{equation}
It is easy to check that 
\[
H_\varphi (x) = \sum_{j = 1}^\infty |x|^{2-n-j}\, \varphi_j(x/|x|)
\]
is the solution of (\ref{Hharmonicc}). Let us fix $|x|$. We have
\begin{equation}\label{stimaH}
|H_\varphi (x) | \leq  \sum_{j = 1}^\infty |x|^{2-n-j}\, |\varphi_j(x/|x|)| = |x|^{1-n}\, |\varphi_1(x/|x|)| + \sum_{j = 2}^\infty |x|^{2-n-j}\, |\varphi_j(x/|x|)|
\end{equation}
Let us try to estimate $\|\varphi_j\|_{L^{\infty}(S^{n-1})}$. From (\ref{decomposition}) we have 
\[
\|\varphi_j\|_{W^{2k,2}(S^{n-1})} \leq c\, j^k\, (n-2+j)^k\, \|\varphi_j\|_{L^{2}(S^{n-1})}
\]
and by the Sobolev embedding theorem we have that $W^{2k,2}(S^{n-1}) \subseteq L^{\infty}(S^{n-1})$ when $4k > n-1$.
We conclude that there exists a positive number $P(n)$ depending only on the dimension $n$ such that 
\[
\|\varphi_j\|_{L^{\infty}(S^{n-1})} \leq c\, j^{P(n)}\, \|\varphi_j\|_{L^{2}(S^{n-1})}
\]
Moreover
\[
\|\varphi_j\|^2_{L^{2}(S^{n-1})} \leq \|\varphi\|^2_{L^{2}(S^{n-1})} \leq \textnormal{ Vol}_{\mathring g}(S^{n-1})\, \|\varphi\|^2_{L^{\infty}(S^{n-1})}
\]
and we can conclude that there exists a constant $c$ such that 
\[
\|\varphi_j\|_{L^{\infty}(S^{n-1})} \leq c\, j^{P(n)}\, \|\varphi\|_{L^{\infty}(S^{n-1})}
\]
From (\ref{stimaH}) we get
\[
|H_\varphi (x) | \leq c\, |x|^{1-n}\, \|\varphi\|_{L^{\infty}(S^{n-1})}\, \left(1 + \sum_{j = 2}^\infty |x|^{1-j}\, j^{P(n)} \right)
\]
It is easy to check that for $|x| \geq 2$
\[
\sum_{j = 2}^\infty |x|^{1-j}\, j^{P(n)} \leq \infty
\]
and this allows us to conclude that for $|x| \geq 2$ there exists a constant $c$ such that 
\begin{equation}\label{firstest}
|H_\varphi (x) | \leq c\, |x|^{1-n}\, \|\varphi\|_{L^{\infty}(S^{n-1})}
\end{equation}
By the maximum principle this inequality is valid also for $1 \leq |x| \leq 2$. Standard elliptic estimates apply to give also 
\begin{equation}\label{secondest}
|\nabla H_{\varphi}(x)| \leq c\, |x|^{-n}\, \| \varphi \|_{L^{\infty}(S^{n-1})}
\end{equation}
Finally, (\ref{firstest}) and (\ref{secondest}) give the following estimate
\[
\|H_{\varphi}(x)\|_{\mathcal{C}^{2,\alpha}_{1-n}(\mathbb{R}^{n} \setminus \mathring B_1)} \leq c\, \| \varphi \|_{\mathcal{C}^{2,\alpha}(S^{n-1})}
	\]
for some constant $c$. From (\ref{firstest}) it is clear that 
\[
\lim_{|x| \rightarrow + \infty} H_{\varphi}(x) = 0.
\]
This completes the proof of the Lemma. \hfill $\Box$
\medskip

Let us define a continuous extension of $H_{\varphi}$ to $\mathbb{R}^n$ in this way :
\begin{equation}
\tilde H_{\varphi} (x) = 
\left\{
\begin{array}{lll}
	0 & \textnormal{for} & |x| \leq \frac{1}{2}\\[3mm]
	\left(2|x| - 1\right) H_{\varphi}\left(\frac{x}{|x|}\right) & \textnormal{for} & \frac{1}{2} \leq |x| \leq 1\\[3mm]
	H_{\varphi}(x) & \textnormal{for} & \mathbb{R}^{n} \setminus \mathring B_1
\end{array}
\right.
\end{equation}
and let us denote
\[
H_{\varphi, \epsilon} = H_{\varphi}\left(\frac{x}{\epsilon}\right)
\]
and
\[
\tilde H_{\varphi, \epsilon} = \tilde H_{\varphi}\left(\frac{x}{\epsilon}\right).
\]

\medskip 

Let $\chi$ be a cutoff function identically equal to 1 for $|x| \leq R_0 << 1$ ($R_0$ is chosen in such a way that $B_{R_0}^g(p)$ belongs to the $x$-coordinate neighborhood of $p$) and identically equal to 0 in $M \setminus B_{2R_0}^g(p)$. 

\medskip 

The main result of this section is the following :

\begin{Proposition}\label{eigenfunction}
Let us suppose $n \geq 3$ and $\nu \in (2-n, \min\{4-n,0\})$. For all $\epsilon$ small enough there exist $(\Lambda_\epsilon,\varphi_\epsilon,w_\epsilon)$ in a neighborhood of $(0,0,0)$ in $\mathbb{R} \times \mathcal{C}^{2,\alpha}_m(S^{n-1}) \times \mathcal{C}_{\nu}^{2,\alpha}(M \setminus \{p\}))$ such that the function
\begin{equation}\label{phi1}
\phi_\epsilon = \phi_0 - \epsilon^{n-2}\,(\phi_0(p) + \Lambda_\epsilon)\,\Gamma_p + w_\epsilon + \chi\, \tilde H_{\varphi_\epsilon, \epsilon}
\end{equation}
(considered in $M \setminus B_\epsilon^g(p)$), is a positive solution of (\ref{foutball}) where
\begin{equation}
\lambda = \lambda_0 + \epsilon^{n-2}\,\mu
\end{equation}
with 
\begin{equation}\label{mu}
\mu = c_n\, \phi_0(p)^2 + \mathcal{O}(\epsilon).
\end{equation}
Moreover the following estimations hold :
\begin{itemize}
\item If $\phi_0$ is not a constant function (CASE 1) then there exists a positive constant $c$ such that
	\[
|\Lambda_\epsilon| + \left\|\varphi_\epsilon\right\|_{L^{\infty}(S^{n-1})} \leq c\, \epsilon\, 
	\qquad \textnormal{and} \qquad \| w_\epsilon \|_{\mathcal{C}_{\nu}^{2,\alpha}(M \setminus \{p\})} \leq c\, (\epsilon^{2n-4} + \epsilon^n + \epsilon^{3-\nu})
	\]
	\item If $\phi_0$ is a constant function (CASE 2) then there exists a positive constant $c$ such that
	\[
\begin{array}{ll}
|\Lambda_\epsilon| + \left\|\varphi_\epsilon\right\|_{L^{\infty}(S^{n-1})} \leq c \,\epsilon & \textnormal{if}\,\,n = 3\\[3mm]
|\Lambda_\epsilon| + \left\|\varphi_\epsilon\right\|_{L^{\infty}(S^{n-1})} \leq c \,\epsilon^\beta\,\,\, \forall \beta < 2 & \textnormal{if}\,\,n = 4\\[3mm]
|\Lambda_\epsilon| + \left\|\varphi_\epsilon\right\|_{L^{\infty}(S^{n-1})} \leq c \,\epsilon^2 & \textnormal{if}\,\,n \geq 5
\end{array}
\]
	and 
	\[
\begin{array}{ll}
\|w_\epsilon\|_{\mathcal{C}^{2,\alpha}_\nu(M \setminus \{p\})} \leq c\, \epsilon^2 & \textnormal{if}\,\,n = 3\\[3mm]
\|w_\epsilon\|_{\mathcal{C}^{2,\alpha}_\nu(M \setminus \{p\})} \leq c\, \left(\epsilon^{4} + \epsilon^{\beta-\nu}\right)\,\,\, \forall \beta < 4 & \textnormal{if}\,\,n = 4\\[3mm]
\|w_\epsilon\|_{\mathcal{C}^{2,\alpha}_\nu(M \setminus \{p\})} \leq c\, \left(\epsilon^{2n-4} + \epsilon^{1+n} + \epsilon^{4-\nu}\right) & \textnormal{if}\,\,n \geq 5
\end{array}
\]
	
\end{itemize}

\end{Proposition}

{\bf Proof.} First we prove that 
\begin{equation}\label{n-2}
\lambda - \lambda_0 = \mathcal{O}(\epsilon^{n-2}).
\end{equation}
By definition
\begin{equation}\label{lH}
\lambda = \min_{u \in H_0^1(M \setminus B_{\epsilon}^g(p))} \frac{\displaystyle \int_{M \setminus B_{\epsilon}^g(p)} |\nabla^g u|^2\, \mbox{dvol}_{g}}{\displaystyle \int_{M \setminus B_{\epsilon}^g(p)} u^2\, \mbox{dvol}_{g}}
\end{equation}
Let us consider a sequence of functions $u_j \in H_0^1(M \setminus B_{\epsilon}^g(p))$ converging to the function
\[
u_* (x) = 
\left\{
\begin{array}{lll}
	\displaystyle \left(\frac{|x|}{\epsilon} - 1\right) \phi_0\left(\frac{2\, \epsilon\, x}{|x|}\right) & \textnormal{in} & B_{2\,\epsilon}^g(p) \setminus B_\epsilon^g(p)\\[3mm]
	\displaystyle \phi_0 & \textnormal{in} & M \setminus B_{2\,\epsilon}^g(p)
\end{array}
\right.
\]
It is easy to check that
\[
\int_{M \setminus B_{\epsilon}^g(p)} u_*^2\, \mbox{dvol}_{g}= \int_{M} \phi_0^2\, \mbox{dvol}_{g}+ \mathcal{O}(\epsilon^n)
\]
while
\[
\int_{M \setminus B_{\epsilon}^g(p)} |\nabla u_*|^2\, \mbox{dvol}_{g}= \int_{M} |\nabla \phi_0|^2\, \mbox{dvol}_{g} + \mathcal{O}(\epsilon^{n-2})
\]
then, from the last two relations and (\ref{lH}), it follows (\ref{n-2}). This allows us to search $\lambda$ in the form 
\begin{equation}\label{lambdaa}
\lambda = \lambda_0 + \epsilon^{n-2}\, \mu
\end{equation}
where $\mu = o(1)$.

\medskip

Let us chose $\phi_\epsilon$ in the form
\[
\phi_\epsilon = \phi_0 - \epsilon^{n-2}\,(\phi_0(p) + \Lambda)\,\Gamma_p + w + \chi\, \tilde H_{\varphi, \epsilon}
\]
for some $(\Lambda,\varphi,w) \in \mathbb{R} \times \mathcal{C}^{2,\alpha}_m(S^{n-1}) \times \mathcal{C}^{2,\alpha}(M \setminus \{p\})$. Then $\phi_\epsilon$ satisfy the first equation of (\ref{foutball}) over $M \setminus B^g_\epsilon(p)$, with $\lambda$ as in (\ref{lambdaa}), if and only if :
\begin{equation}\label{ww1}
\begin{array}{lcl}
\Big( \Delta_g + \lambda_0 + \epsilon^{n-2}\, \mu \Big)\, w + \epsilon^{n-2}\, \Big[ \mu - c_n\, \phi_0(p)\, (\phi_0(p) + \Lambda)\Big]\, \phi_0 + H_{\varphi, \epsilon}\, \Delta_g \chi + & &\\[3mm]
\qquad + \chi\, \Delta_g H_{\varphi,\epsilon} + 2\, \nabla^g H_{\varphi, \epsilon}\, \nabla^g \chi - \epsilon^{2n-4}\, \mu \, (\phi_0(p) + \Lambda)\, \Gamma_p + (\lambda_0 + \epsilon^{n-2}\, \mu)\, \chi\,H_{\varphi, \epsilon} & = & 0
\end{array}
\end{equation}
over $M \setminus B^g_\epsilon(p)$.
This equation can be considered over $M \setminus \{p\}$ if we take $\tilde H_{\varphi,\epsilon}$ instead of $H_{\varphi,\epsilon}$, and a continuous extension $\widetilde{\Delta_g H}_{\varphi,\epsilon}$ of the function $\Delta_g H_{\varphi,\epsilon}$ (similarly to the continuous extension $\tilde H_{\varphi,\epsilon}$ of the function $H_{\varphi,\epsilon}$). Remark that the term $\nabla^g H_{\varphi, \epsilon}\, \nabla^g \chi$ is 0 in a neighborhood of $\partial B^g_\epsilon(p)$ then can be extended to 0 in $B^g_\epsilon(p)$.

\medskip

We need the following :
\begin{Lemma}\label{espacepoids}
Let $n \geq 3$. The operator
\[
\big(\Delta_g + \lambda_0 +\epsilon^{n-2}\,\mu\big): \mathcal{C}^{2,\alpha}_{\nu, \bot,0}(M \setminus \{p\}) \longrightarrow \mathcal{C}^{0,\alpha}_{\nu-2,\bot}(M \setminus \{p\}),
\]
where the subscript $\bot$ is meant to point out that functions are $L^2$-orthogonal to $\phi_0$ and the subscript $0$ is meant to point out that functions satisfy the 0 Dirichlet (in the CASE 1) or 0 Neumann (in the CASE 2) boundary condition on $\partial M$ if $\partial M \neq \emptyset$, 
 is an isomorphism for $\nu \in (2-n,0)$ and $\epsilon$ small enough.
\end{Lemma}

{\bf Proof.} Let us suppose that $\nu \in (2-n,0)$ and $n \geq 3$. In \cite{notepacard} is proved that for all $f \in \mathcal{C}^{0,\alpha}_{\nu-2}(\mathring B_1 \setminus \{0\})$ there exists a unique solution $u \in \mathcal{C}^{2,\alpha}_{\nu}(\mathring B_1 \setminus \{0\})$ of 
\begin{equation}
\left\{
\begin{array}{rcccl}
	\Delta_{\mathring g} \, u & = & f & \textnormal{in} & \mathring B_1 \setminus \{0\} \\[3mm]
	u & = & 0 & \textnormal{on} & \partial \mathring B_1
\end{array}.
\right.
\end{equation}
 
\medskip 

Now let us consider the operator $\Delta_g + \lambda_0$. Let $R_0$ small enough, take the normal geodesic coordinates in $B^g_{R_0}(p)$, and let $f \in \mathcal{C}^{0,\alpha}_{\nu-2}(M \setminus \{p\})$. Considering the dilated metric $R_0^{-2}\, g$, the parameterization of $B^{g}_{R_0}(p)$ given by
\[
Y ( y) : =\mbox{Exp}_p^{g} \left( R_0\, \sum_i y^i \, E_i \right)
\]
and the ball $\mathring B_1$ endowed with the metric $\check g = Y^{*}(R_0^{-2}\, g)$, the problem
\[
\left\{
\begin{array}{rclcl}
	(\Delta_{g} \, + \lambda_0) \, u & = & f & \textnormal{in} & B^g_{R_0} \setminus \{p\} \\[3mm]
	u & = & 0 & \textnormal{on} & \partial B^g_{R_0} 
\end{array}
\right.
\]
is equivalent, in terms of existence and unicity of solutions, to the problem 
\[
\left\{
\begin{array}{rclcl}
	(\Delta_{\check g} \, + R_0^2\, \lambda_0) \, u & = & Y^* f & \textnormal{in} & \mathring B_1 \setminus \{0\} \\[3mm]
	u & = & 0 & \textnormal{on} & \partial \mathring B_1 
\end{array}.
\right.
\]
Considering that the difference between the coefficients of the metric $\check g$ and the metric $\mathring g$ can be estimate by a constant times $R_0^2$, the operator $\Delta_{\check g} + R_0^2\, \lambda_0$ is a small perturbation of the operator $\Delta_{\mathring g}$ when $R_0$ is small. By the previous claim, we conclude that there exists a positive $R_0$ (small enough) such that, when $\nu \in (2-n,0)$ and $n \geq 3$, for all $f \in \mathcal{C}^{0,\alpha}_{\nu-2}(M \setminus \{p\})$ there exists a unique solution $u \in \mathcal{C}^{2,\alpha}_{\nu}(B^g_{R_0} \setminus \{p\})$ of 
\[
\left\{
\begin{array}{rclcl}
	(\Delta_{g} \, + \lambda_0) \, u & = & f & \textnormal{in} & B^g_{R_0} \setminus \{p\} \\[3mm]
	u & = & 0 & \textnormal{on} & \partial B^g_{R_0} 
\end{array}
\right.
\]

\medskip

Let now consider the solution of 
\begin{equation}\label{vvv}
(\Delta_{g} \, + \lambda_0) \, v = f - (\Delta_g + \lambda_0)\,(\tilde \chi\, u)
\end{equation}
with 0 Dirichlet boundary condition at $\partial M$, where $\tilde \chi$ is a cut-off function equal to 1 for $|x| \leq R_0/2$ and equal to 0 for $|x| \geq R_0$. We remark that this equation is well defined in $M$, because the singularity of $f$ at $p$ is balanced by $(\Delta_g + \lambda_0)\,(\tilde \chi\, u)$. Moreover the right hand side term is orthogonal to $\phi_0$ if $f$ has such a property. Hence, there exists a solution $v \in \mathcal{C}^{2,\alpha}_{\bot,0}(M)$ to (\ref{vvv}), and we have that 
\[
	(\Delta_{g} \, + \lambda_0) \, (\tilde \chi\, u + v)  = f 
\]
in $M \setminus \{p\}$, with 0 Dirichlet condition at $\partial M$. Obviously $w = \tilde \chi\, u + v  \in \mathcal{C}^{2,\alpha}_{\nu,\bot}(M \setminus \{p\})$. We proved then that for $\nu \in (2-n,0)$ and $n \geq 3$ and for all $f \in \mathcal{C}^{0,\alpha}_{\nu-2,\bot}(M \setminus \{p\})$ there exists a unique solution $w \in \mathcal{C}^{2,\alpha}_{\nu,\bot}(M \setminus \{p\})$ of 
\[
	(\Delta_{g} \, + \lambda_0) \, w  = f 
\]
in $M \setminus \{p\}$ with 0 Dirichlet at $\partial M$. This result is still true for the operator $\Delta_g + \lambda_0 + \epsilon^{n-2}\, \nu$ when $\epsilon$ is small enough, because such an operator is a small perturbation of the operator $\Delta_{g} \, + \lambda_0$. The raisoning done do not change if we consider the 0 Neumann boundary condition on $\partial M$ instead of the 0 Dirichlet boundary condition. This completes the proof of the Lemma. \hfill $\Box$

\medskip

To semplify the notation let us define
\[
\begin{array}{lcl}
A &:=& \epsilon^{n-2}\, \Big[ \mu - c_n\, \phi_0(p)\, (\phi_0(p) + \Lambda)\Big]\, \phi_0\\[3mm]
B &:=& \tilde H_{\varphi, \epsilon}\, \Delta_g \chi + \chi\, \widetilde{\Delta_g  H}_{\varphi,\epsilon}
+ 2\, \nabla^g  H_{\varphi, \epsilon}\, \nabla^g \chi\\[3mm]
C &:=& - \epsilon^{2n-4}\, \mu \, (\phi_0(p) + \Lambda)\, \Gamma_p\\[3mm]
D &:=& (\lambda_0 + \epsilon^{n-2}\, \mu)\, \chi\, \tilde H_{\varphi, \epsilon}
\end{array}
\]
We remark that $\Gamma_p \in \mathcal{C}^{0,\alpha}_{\nu-2}(M \setminus \{p\})$ if $\nu < 4-n$.
Equation (\ref{ww1}), extended to $M \backslash \{p\}$, becomes
\[
( \Delta_g + \lambda_0 + \epsilon^{n-2}\, \mu )\, w = - (A + B +C+D)
\]
By the last result, if we chose $\mu$ in order to verify 
\begin{equation}\label{ABCD1}
\int_M (A + B +C+D)\, \phi_0 = 0
\end{equation}
there exists a solution $w(\epsilon,\Lambda,\varphi) \in \mathcal{C}^{2,\alpha}_{\nu,\bot,0}(M \setminus \{p\})$ to equation (\ref{ww1}) for 
\[
\nu \in \left(2-n,\min\{0,4-n\}\right),
\]
for all $\Lambda \in \mathbb{R}$, for all $\varphi \in  \mathcal{C}^{2,\alpha}_m(S^{n-1})$, and for all $\epsilon$ small enough, and then
\[
\phi_\epsilon = \phi_0 + \epsilon^{n-2}\, (\phi_0(p) + \Lambda)\, \Gamma_p + w(\epsilon, \Lambda, \varphi) + \chi\, H_{\varphi, \epsilon}
\]
with $w(\epsilon, \Lambda, \varphi)$ restricted to $M \backslash B^g_\epsilon(p)$,
satisfies the first equation of (\ref{foutball}). From (\ref{ABCD1}) we get
\[
\mu = \frac{\displaystyle \epsilon^{n-2}\, c_n\, \phi_0(p)\, (\phi_0(p) + \Lambda) - \int_M B\, \phi_0 - \lambda_0\, \int_M \chi\, \tilde H_{\varphi, \epsilon}\, \phi_0}{\displaystyle \epsilon^{n-2} \left( 1 + \int_M \chi\, \tilde H_{\varphi, \epsilon}\, \phi_0\right)}
\]
It is easy to check that 
\[
\int_M B\, \phi_0 \leq c\, \epsilon^{n-1}\, \|\varphi\|_{L^{\infty}(S^{n-1})}
\]
and 
\[
\int_M \chi\, \tilde H_{\varphi, \epsilon}\, \phi_0 \leq c\, \epsilon^{n-1}\, \|\varphi\|_{L^{\infty}(S^{n-1})}
\]
from which it follows the expansion of $\mu$:
\begin{equation}\label{muexpansion}
\mu = c_n\, \phi_0(p)\, (\phi_0(p) + \Lambda) + \mathcal{O}(\epsilon)\, \|\varphi\|_{L^\infty(S^{n-1})}.
\end{equation}

\medskip

We want now to give some estimations on the function $w$. By the previous results and Lemma \ref{lemmaharmonic} we have the following estimations :
\medskip

\begin{itemize}
	\item $\| A \|_{\mathcal{C}^{0,\alpha}_{\nu-2}(M \setminus \{p\})} \leq c\, \epsilon^{n-1}\,\|\varphi\|_{L^{\infty}(S^{n-1})}$
	\medskip
	
	\item $\| B \|_{\mathcal{C}^{0,\alpha}_{\nu-2}(M \setminus \{p\})} \leq c\,(\epsilon^{n-1}+ \epsilon^{2-\nu})\, \|\varphi\|_{L^{\infty}(S^{n-1})}$ 
	\medskip
	\item $\| C \|_{\mathcal{C}^{0,\alpha}_{\nu-2}(M \setminus \{p\})} \leq c\, \epsilon^{2n-4}$
	\medskip
	\item $\| D \|_{\mathcal{C}^{0,\alpha}_{\nu-2}(M \setminus \{p\})} \leq c\, \left(\epsilon^{2-\nu} + \epsilon^{n-1}\right) \, \|\varphi\|_{L^{\infty}(S^{n-1})}$
	\medskip
\end{itemize}
In particular we get
\[
\| A + B+ C + D \|_{\mathcal{C}^{0,\alpha}_{\nu-2}(M \setminus \{p\})} \leq c\, \left( 
\epsilon^{2n-4} + \epsilon^{n-1}\, \|\varphi\|_{L^{\infty}(S^{n-1})} + \epsilon^{2-\nu}\, \|\varphi\|_{L^{\infty}(S^{n-1})}\right)
\]
This give us an estimation for the function $w$ that we found before:
\[
\| w \|_{\mathcal{C}^{2,\alpha}_{\nu}(M \setminus \{p\})} \leq c\, \left(
 \epsilon^{2n-4} + \epsilon^{n-1}\, \|\varphi\|_{L^{\infty}(S^{n-1})} + \epsilon^{2-\nu}\, \|\varphi\|_{L^{\infty}(S^{n-1})}\right).
\]



\medskip 

We have proved the following :

\medskip

\textbf{First intermediate result. }\textit{Let $\nu \in (2-n,4-n)$. For all $\Lambda \in \mathbb{R}$, for all $\varphi \in \mathcal{C}^{2,\alpha}_{m}(S^{n-1})$, for all $\epsilon$ small enough, there exists a function $w(\epsilon, \Lambda, \varphi) \in \mathcal{C}^{2,\alpha}_{\nu,\bot, 0}(M \setminus \{p\})$ such that (\ref{phi1}) is a positive solution of the first equation of (\ref{foutball}). Moreover there exists a positive constant $c$ such that 
\begin{equation}\label{estimationdew}
\| w \|_{\mathcal{C}^{2,\alpha}_{\nu}(M \setminus \{p\})} \leq c\, \left( \epsilon^{2n-4} + \epsilon^{n-1}\, \|\varphi\|_{L^{\infty}(S^{n-1})} + \epsilon^{2-\nu}\, \|\varphi\|_{L^{\infty}(S^{n-1})}\right).
\end{equation}}
\medskip

Now we have to make attention to the second equation of (\ref{foutball}). Let us define
\[
N(\epsilon, \Lambda,\varphi) := \Big[ \phi_0(\epsilon\, y) - \epsilon^{n-2}\, (\phi_0(p) + \Lambda)\, \Gamma(\epsilon\, y) + \big(w(\epsilon,\Lambda,\varphi)\big)(\epsilon\, y) + \varphi(y)\Big]_{y \in S^{n-1}}.
\]
We remark that $N$ represents the boundary value of the solution of the first equation of (\ref{foutball}) that we found above, is well defined in a heighborhood of $(0,0,0)$ in $(0,+\infty) \times \mathbb{R} \times \mathcal{C}^{2,\alpha}_m(S^{n-1})$, and takes its values in $\mathcal{C}^{2,\alpha}(S^{n-1})$. It is easy to compute the differential of $N$ with respect to $\Lambda$ and $\varphi$ at $(0,0,0)$ :
\[
\big(\partial_{\Lambda} N (0,0,0)\big)(\tilde \Lambda) = - \tilde \Lambda\\[1mm]
\]
\[
\big(\partial_{\varphi} N (0,0,0)\big)(\tilde \varphi) = \tilde \varphi.
\]
From the estimation of the function $w$ it follows that 
\begin{eqnarray*}
\| w\|_{L^{\infty}(\partial B_{\epsilon}^g(p))} &\leq& \epsilon^{\nu}\, \| w \|_{\mathcal{C}^{2,\alpha}_{\nu}(M \setminus \{p\})}\\
&\leq& c\, \left( \epsilon^{2n-4+\nu}  + \epsilon^{n-1+\nu}\, \|\varphi\|_{L^{\infty}(S^{n-1})} + \epsilon^{2}\, \|\varphi\|_{L^{\infty}(S^{n-1})}\right).
\end{eqnarray*}
Then we can estimate $N(\epsilon, 0,0)$ :
\[
\left\|N(\epsilon, 0,0)\right\|_{L^{\infty}(S^{n-1})} \leq \left\|\phi_0(\epsilon\, x) - \epsilon^{n-2}\, \phi_0(p)\, \Gamma(\epsilon\, x)\right\|_{L^{\infty}(S^{n-1})} + \left\|\big(w(\epsilon,0,0)\big)(\epsilon\, x)\right\|_{L^{\infty}(S^{n-1})}
\]
Here we have again to distinguish the different cases, according on the behaviour of the function $\phi_0$. If $\phi_0$ is not a constant function (CASE 1) we have (using the expansion (\ref{gammaexpansion}) of $\Gamma_p$)
\[
\left\|N(\epsilon, 0,0)\right\|_{L^{\infty}(S^{n-1})} \leq c\, \epsilon . 
\]
The same estimate is obtained if $\phi_0$ is a constant function (CASE 2) and $n=3$. In the CASE 2 and $n=4$ we get
\[
\left\|N(\epsilon, 0,0)\right\|_{L^{\infty}(S^{n-1})} \leq c\, \epsilon^\beta 
\]
$\forall \beta < 2$ and when $n \geq 5$ :
\[
\left\|N(\epsilon, 0,0)\right\|_{L^{\infty}(S^{n-1})} \leq c\, \epsilon^2 
\]
The implicit function theorem applies to give the :

\medskip

\textbf{Second intermediate result. }\textit{Let $\nu \in (2-n,\min\{4-n,0\})$, let $w$ be the function found in the first intermediate result, and let $\epsilon$ be small enough. Then there exist $(\Lambda_\epsilon,\varphi_\epsilon)$ in a neighborhood of $(0,0)$ in $\mathbb{R} \times \mathcal{C}^{2,\alpha}_m(S^{n-1})$ such that $N(\epsilon, \Lambda_\epsilon, \varphi_\epsilon) = 0$ (i.e. (\ref{phi1}) is a positive solution of (\ref{foutball})). Moreover the following estimations hold :
\begin{itemize}
\item If $\phi_0$ is not a constant function (CASE 1) then 
	\[
	|\Lambda_\epsilon| + \left\|\varphi_\epsilon \right\|_{L^{\infty}(S^{n-1})} \leq c \,\epsilon 
	\]
	\item If $\phi_0$ is a constant function (CASE 2) then 
	\[
\begin{array}{ll}
|\Lambda_\epsilon| + \left\|\varphi_\epsilon\right\|_{L^{\infty}(S^{n-1})} \leq c \,\epsilon & \textnormal{if}\,\,n = 3\\[3mm]
|\Lambda_\epsilon| + \left\|\varphi_\epsilon\right\|_{L^{\infty}(S^{n-1})} \leq c \,\epsilon^\beta\,\,\, \forall \beta < 2 & \textnormal{if}\,\,n = 4\\[3mm]
|\Lambda_\epsilon| + \left\|\varphi_\epsilon\right\|_{L^{\infty}(S^{n-1})} \leq c \,\epsilon^2 & \textnormal{if}\,\,n \geq 5
\end{array}
\]
\end{itemize}
}

\medskip 

Putting together the first and the second intermediate result, we get the following existence result: for all $\epsilon$ small enough there exist $(\Lambda_\epsilon,\varphi_\epsilon,w_\epsilon)$ in a neighborhood of $(0,0,0)$ in $\mathbb{R} \times \mathcal{C}^{2,\alpha}_m(S^{n-1}) \times \mathcal{C}^{2,\alpha}(M \setminus \{p\}))$ such that the function
\[
\phi_\epsilon = \phi_0 - \epsilon^{n-2}\,(\phi_0(p) + \Lambda_\epsilon)\,\Gamma_p + w_\epsilon + \chi\, H_{\varphi_\epsilon, \epsilon}
\]
considered in $M \setminus B_\epsilon^g(p)$, is a positive solution of (\ref{foutball}) where
\begin{equation}
\lambda = \lambda_0 + \epsilon^{n-2}\, \mu
\end{equation}
and from (\ref{muexpansion}) we obtain 
\begin{equation}\label{mu}
\mu = c_n\, \phi_0(p)^2 + \mathcal{O}(\epsilon).
\end{equation}
Moreover :
\begin{itemize}
\item If $\phi_0$ is not a constant function (CASE 1) then there exists a positive constant $c$ such that
	\[
	|\Lambda_\epsilon| + \left\|\varphi_\epsilon\right\|_{L^{\infty}(S^{n-1})} \leq c\, \epsilon
	\]
and from (\ref{estimationdew})
	\[
\| w_\epsilon \|_{\mathcal{C}^{2,\alpha}_\nu(M \setminus \{p\}))}\leq c\,\left( \epsilon^{2n-4} + \epsilon^n + \epsilon^{3-\nu}\right)
	\]
	\item If $\phi_0$ is a constant function (CASE 2) then there exists a positive constant $c$ such that
	\[
\begin{array}{ll}
|\Lambda_\epsilon| + \left\|\varphi_\epsilon\right\|_{L^{\infty}(S^{n-1})} \leq c \,\epsilon & \textnormal{if}\,\,n = 3\\[3mm]
|\Lambda_\epsilon| + \left\|\varphi_\epsilon\right\|_{L^{\infty}(S^{n-1})} \leq c \,\epsilon^\beta\,\,\, \forall \beta < 2 & \textnormal{if}\,\,n = 4\\[3mm]
|\Lambda_\epsilon| + \left\|\varphi_\epsilon\right\|_{L^{\infty}(S^{n-1})} \leq c \,\epsilon^2 & \textnormal{if}\,\,n \geq 5
\end{array}
\]
	and from (\ref{estimationdew})
	\[
\begin{array}{ll}
\|w_\epsilon\|_{\mathcal{C}^{2,\alpha}_\nu(M \setminus \{p\})} \leq c\, \epsilon^2 & \textnormal{if}\,\,n = 3\\[3mm]
\|w_\epsilon\|_{\mathcal{C}^{2,\alpha}_\nu(M \setminus \{p\})} \leq c\, \left(\epsilon^{4} + \epsilon^{\beta-\nu}\right)\,\,\, \forall \beta < 4 & \textnormal{if}\,\,n = 4\\[3mm]
\|w_\epsilon\|_{\mathcal{C}^{2,\alpha}_\nu(M \setminus \{p\})} \leq c\, \left(\epsilon^{2n-4} + \epsilon^{1+n} + \epsilon^{4-\nu}\right) & \textnormal{if}\,\,n \geq 5
\end{array}
\]
\end{itemize}

This completes the proof of the result. 
\hfill $\Box$

\medskip

For the case $n=2$ we can generalise the previuos proposition, obtaining the:

\begin{Proposition}\label{eigenfunction2}
Let we suppose $n =2$ and $\nu \in (0,1)$. For all $\epsilon$ small enough there exist $(\Lambda_\epsilon,\varphi_\epsilon,w_\epsilon)$ in a neighborhood of $(0,0,0)$ in $\mathbb{R} \times \mathcal{C}^{2,\alpha}_m(S^{n-1}) \times \Big( \tilde \chi\, \mathbb{R} \oplus \mathcal{C}^{2,\alpha}_\nu(M \setminus \{p\})\Big)$, where $\tilde \chi$ is some cut-off funtion equal to 1 in a neighborhood of the origin, 
such that the function
\begin{equation}\label{phi}
\phi_\epsilon = \phi_0 - (\phi_0(p) + \Lambda_\epsilon)\,\Gamma_p + w_\epsilon + \chi\, \tilde H_{\varphi_\epsilon, \epsilon}
\end{equation}
considered in $M \setminus B_\epsilon^g(p)$, is a positive solution of (\ref{foutball}) where
\begin{equation}
\lambda = \lambda_0 + (\log \epsilon)^{-1}\,\mu
\end{equation}
with 
\begin{equation}\label{mu}
\mu = c_n\, \phi_0(p)^2 + \mathcal{O}(\epsilon).
\end{equation}
Moreover the following estimations hold : there exists a positive constant $c$ such that
	\[
|\Lambda_\epsilon| + \left\|\varphi_\epsilon\right\|_{L^{\infty}(S^{n-1})} \leq c\, \epsilon 
	\qquad \textnormal{and} \qquad \| w_\epsilon \|_{\tilde \chi\, \mathbb{R}\, \oplus\, \mathcal{C}^{2,\alpha}_\nu(M \setminus \{p\})} \leq c\, \epsilon^2\, \log\epsilon
	\]
	

\end{Proposition}

{\bf Proof.} We will follow the proof of the previous proposition, adapting it to the case of dimension 2. First we prove that 
\begin{equation}\label{log}
\lambda - \lambda_0 \leq c\, \log(\epsilon)^{-1}.
\end{equation}
Let we consider a sequence of functions $u_j \in H_0^1(M \setminus B_{\epsilon}^g(p))$ converging to the function
\[
u_* (x) = 
\left\{
\begin{array}{lll}
	\displaystyle K\, \log \frac{|x|}{\epsilon} & \textnormal{in} & B_{R_0}^g(p) \setminus B_\epsilon^g(p)\\[3mm]
	\displaystyle \phi_0 & \textnormal{in} & M \setminus B_{R_0}^g(p)
\end{array}
\right.
\]
where $K$ is a constant chosen in order to make the function $u_*$ continuous.
It is easy to check that for $\epsilon$ small enough
\[
\int_{M \setminus B_{\epsilon}^g(p)} u_*^2\, \mbox{dvol}_{g}= \leq c\, \log(\epsilon)^2
\]
while
\[
\int_{M \setminus B_{\epsilon}^g(p)} |\nabla u_*|^2\, \mbox{dvol}_{g} \leq c\, \log(\epsilon)
\]
then, from the last two relations and (\ref{lH}), it follows (\ref{log}). This allows us to search $\lambda$ in the form 
\begin{equation}\label{lambda}
\lambda = \lambda_0 + \log (\epsilon)^{-1}\, \mu
\end{equation}
where $\mu = o(1)$.
\medskip

Let we chose $\phi_\epsilon$ in the form
\[
\phi_\epsilon = \phi_0 - (\log \epsilon)^{-1}\,(\phi_0(p) + \Lambda)\,\Gamma_p + w + \chi\, \tilde H_{\varphi, \epsilon}
\]
for some $(\Lambda,\varphi,w) \in \mathbb{R} \times \mathcal{C}^{2,\alpha}_m(S^{n-1}) \times \mathcal{C}^{2,\alpha}(M \setminus \{p\})$. Then $\phi_\epsilon$ satisfy the first equation of (\ref{foutball}), with $\lambda$ as in (\ref{lambda}), if and only if the quantity
\begin{equation}\label{w}
\begin{array}{lcl}
\Big( \Delta_g + \lambda_0 + (\log \epsilon)^{-1}\, \mu \Big)\, w + (\log \epsilon)^{-1}\, \Big[ \mu - c_n\, \phi_0(p)\, (\phi_0(p) + \Lambda)\Big]\, \phi_0 + \tilde H_{\varphi, \epsilon}\, \Delta_g \chi + & &\\[3mm]
\qquad + \chi\, \Delta_g \tilde H_{\varphi,\epsilon} + 2\, \nabla^g \tilde H_{\varphi, \epsilon}\, \nabla^g \chi - (\log \epsilon)^{-2}\, \mu \, (\phi_0(p) + \Lambda)\, \Gamma_p + (\lambda_0 + (\log \epsilon)^{-1}\, \mu)\, \chi\, \tilde H_{\varphi, \epsilon} 
\end{array}
\end{equation}
is identically equal to 0 over $M \setminus B_\epsilon^g(p)$. This equation can be consider, after opportune extension of functions like in equation (\ref{ww1}), over $M \setminus \{p\}$.

\medskip
The operator
\[
\big(\Delta_g + \lambda_0 \big): \tilde \chi\,\mathbb{R} \oplus \mathcal{C}^{2,\alpha}_{\nu, \bot,0}(M \setminus \{p\}) \longrightarrow \mathcal{C}^{0,\alpha}_{\nu-2,\bot}(M \setminus \{p\}),
\]
where the subscript $\bot$ is meant to point out that functions are $L^2$-orthogonal to $\phi_0$ and the subscript $0$ is meant to point out that functions satisfy the 0 Dirichlet (in the CASE 1) or 0 Neumann (in the CASE 2) boundary condition on $\partial M$ if $\partial M \neq \emptyset$, and $\tilde \chi$ is some cut-off funtion equal to 1 in a neighborhood of the origin, is an isomorphism for $\nu \in (0,1)$. The same result holds for the operator 
\[
\big(\Delta_g + \lambda_0 +(\log \epsilon)^{-1}\,\mu\big): \tilde \chi\, \mathbb{R} \oplus \mathcal{C}^{2,\alpha}_{\nu, \bot,0}(M \setminus \{p\}) \longrightarrow \mathcal{C}^{0,\alpha}_{\nu-2,\bot}(M \setminus \{p\}),
\]
when $\epsilon$ is small enough. The proof of this fact can be obtained by comparison of the proof of the Lemma \ref{espacepoids} and the analysis of the case of dimension 2 in \cite{notepacard}.

\medskip 

To semplify the notation let we define
\[
\begin{array}{lcl}
A &:=& (\log \epsilon)^{-1}\, \Big[ \mu - c_n\, \phi_0(p)\, (\phi_0(p) + \Lambda)\Big]\, \phi_0\\[3mm]
B &:=& \tilde H_{\varphi, \epsilon}\, \Delta_g \chi + \chi\, \Delta_g \tilde H_{\varphi,\epsilon}
+ 2\, \nabla^g \tilde H_{\varphi, \epsilon}\, \nabla^g \chi\\[3mm]
C &:=& - (\log \epsilon)^{-2}\, \mu \, (\phi_0(p) + \Lambda)\, \Gamma_p\\[3mm]
D &:=& (\lambda_0 + (\log \epsilon)^{-1}\, \mu)\, \chi\, \tilde H_{\varphi, \epsilon}
\end{array}
\]
We remark that $\Gamma \in \mathcal{C}^{0,\alpha}_{\nu-2}(M \setminus \{p\})$ when $\nu \in (0,1)$.
Equation (\ref{w}) becomes
\[
( \Delta_g + \lambda_0 + (\log \epsilon)^{-1}\, \mu )\, w = - (A + B +C+D)
\]
By the last result, if we chose $\mu$ in order to verify
\begin{equation}\label{ABCD}
\int_M (A + B +C+D)\, \phi_0 = 0
\end{equation}
there exists a solution $w(\epsilon,\Lambda,\varphi) = w^{(1)} + w^{(2)} \in \tilde \chi\, \mathbb{R} \oplus \mathcal{C}^{2,\alpha}_{\nu,\bot,0}(M \setminus \{p\})$ to equation (\ref{w}) for $\nu \in (0,1)$, for all $\Lambda \in \mathbb{R}$, for all $\varphi \in  \mathcal{C}^{2,\alpha}_m(S^{n-1})$, and for all $\epsilon$ small enough, and then
\[
\phi_\epsilon = \phi_0 - (\phi_0(p) + \Lambda)\, \Gamma_p + w(\epsilon, \Lambda, \varphi) + \chi\, H_{\varphi, \epsilon}
\]
satisfy the first equation of (\ref{foutball}). From (\ref{ABCD}) we get
\[
\mu = \frac{\displaystyle (\log \epsilon)^{-1}\, c_n\, \phi_0(p)\, (\phi_0(p) + \Lambda) - \int_M B\, \phi_0 - \lambda_0\, \int_M \chi\, \tilde H_{\varphi, \epsilon}\, \phi_0}{\displaystyle (\log \epsilon)^{-1} \left( 1 + \int_M \chi\, \tilde H_{\varphi, \epsilon}\, \phi_0\right)}
\]
It is easy to check that 
\[
\int_M B\, \phi_0 \leq c\, \epsilon\, \|\varphi\|_{L^{\infty}(S^{n-1})}
\]
and 
\[
\int_M \chi\, \tilde H_{\varphi, \epsilon}\, \phi_0 \leq c\, \epsilon\,\|\varphi\|_{L^{\infty}(S^{n-1})}
\]
from which it follows the expansion of $\mu$:
\begin{equation}\label{muexpansion}
\mu = c_n\, \phi_0(p)\, (\phi_0(p) + \Lambda) +  \mathcal{O}(\epsilon\,\log \epsilon)\, \|\varphi\|_{L^\infty(S^{n-1}}.
\end{equation}

\medskip

We want now to give some estimations on the function $w$. By the previous facts and Lemma \ref{lemmaharmonic} we have the following estimations :
\medskip

\begin{itemize}
	\item $\| A \|_{\mathcal{C}^{0,\alpha}_{\nu-2}(M \setminus \{p\})} \leq c\, \epsilon\,\log \epsilon\, \|\varphi\|_{L^{\infty}(S^{n-1})}$
	\medskip
	
	\item $\| B \|_{\mathcal{C}^{0,\alpha}_{\nu-2}(M \setminus \{p\})} \leq c\,\epsilon\, \|\varphi\|_{L^{\infty}(S^{n-1})}$
	\medskip
	\item $\| C \|_{\mathcal{C}^{0,\alpha}_{\nu-2}(M \setminus \{p\})} \leq c\, (\log \epsilon)^{-2}$
	\medskip
	\item $\| D \|_{\mathcal{C}^{0,\alpha}_{\nu-2}(M \setminus \{p\})} \leq c\,(\epsilon+\epsilon^{2-\nu})\, \|\varphi\|_{L^{\infty}(S^{n-1})}$
	\medskip
\end{itemize}
In particular we get
\[
\| A + B+ C + D \|_{\mathcal{C}^{0,\alpha}_{\nu-2}(M \setminus \{p\})} \leq c\, \left(  (\log \epsilon)^{-2}+ \epsilon\,\log \epsilon\, \|\varphi\|_{L^{\infty}(S^{n-1})} \right)
\]
where we used the fact that for $\epsilon$ small enough and $\nu \in (0,1)$ we have $\epsilon^{2-\nu} < \epsilon < \epsilon\, \log\epsilon$.
This give us an estimation on the function $w$ that we found before:
\[
\|w^{(1)}\| + \| w^{(2)} \|_{\mathcal{C}^{2,\alpha}_{\nu}(M \setminus \{p\})} \leq c\, \left( (\log \epsilon)^{-2}+ \epsilon\, \log \epsilon\,  \|\varphi\|_{L^{\infty}(S^{n-1})} \right).
\]

\medskip

We have proved the following :

\medskip

\textbf{First intermediate result. }\textit{Let $\nu \in (0,1)$. For all $\Lambda \in \mathbb{R}$, for all $\varphi \in \mathcal{C}^{2,\alpha}_{m}(S^{n-1})$, for all $\epsilon$ small enough, there exists a function $w(\epsilon, \Lambda, \varphi) = w^{(1)} + w^{(2)} \in \tilde \chi\,\mathbb{R}\oplus\mathcal{C}^{2,\alpha}_{\nu,\bot, 0}(M \setminus \{p\})$ such that (\ref{phi}) is a positive solution of the first equation of (\ref{foutball}). Moreover there exists a positive constant $c$ such that
\[
\|w^{(1)}\| + \| w^{(2)} \|_{\mathcal{C}^{2,\alpha}_{\nu}(M \setminus \{p\})} \leq c\, \left( (\log \epsilon)^{-2}+ \epsilon\, \log \epsilon\, \|\varphi\|_{L^{\infty}(S^{n-1})} \right)
\]}
\medskip

Now we have to make attention to the second equation of (\ref{foutball}). Let we define
\[
N(\epsilon, \Lambda,\varphi) := \Big[ \phi_0(\epsilon\, y) - (\log \epsilon)^{-1}\,(\phi_0(p) + \Lambda)\, \Gamma(\epsilon\, y) + \big(w(\epsilon,\Lambda,\varphi)\big)(\epsilon\, y) + \varphi(y)\Big]_{y \in S^{n-1}}.
\]
We remark that $N$ represents the boundary value of the solution of the first of (\ref{foutball}) we found above, is well defined in a heighborhood of $(0,0,0)$ in $(0,+\infty) \times \mathbb{R} \times \mathcal{C}^{2,\alpha}_m(S^{n-1})$, and takes its values in $\mathcal{C}^{2,\alpha}(S^{n-1})$. The differential of $N$ with respect to $\Lambda$ and $\varphi$ at $(0,0,0)$ is :
\[
\big(\partial_{\Lambda} N (0,0,0)\big)(\tilde \Lambda) = - \tilde \Lambda\\[1mm]
\]
\[
\big(\partial_{\varphi} N (0,0,0)\big)(\tilde \varphi) = \tilde \varphi.
\]
From the estimation of the function $w$ it follows that
\begin{eqnarray*}\label{westimation}
\| w\|_{L^{\infty}(\partial B_{\epsilon}^g(p))} &\leq& 
c\, \left((\log \epsilon)^{-2}+ \epsilon\, \log \epsilon\, \|\varphi\|_{L^{\infty}(S^{n-1})} \right)
\end{eqnarray*}
Then we can estimate $N(\epsilon, 0,0)$ :
\[
\left\|N(\epsilon, 0,0)\right\|_{L^{\infty}(S^{n-1})} \leq \left\|\phi_0(\epsilon\, y) - (\log\epsilon)^{-1}\, \phi_0(p)\, \Gamma(\epsilon\, y)\right\|_{L^{\infty}(S^{n-1})} + \left\|\big(w(\epsilon,0,0)\big)(\epsilon\, y)\right\|_{L^{\infty}(S^{n-1})}
\]
and we have
\[
\left\|N(\epsilon, 0,0)\right\|_{L^{\infty}(S^{n-1})} \leq c\, \epsilon 
\]
The implicit function theorem applies to give the :

\medskip

\textbf{Second intermediate result. }\textit{Let $\nu \in (0,1)$ and let $w$ be the function found in the first intermediate result, and let $\epsilon$ be small enough. Then there exist $(\Lambda_\epsilon,\varphi_\epsilon)$ in a neighborhood of $(0,0)$ is $\mathbb{R} \times \mathcal{C}^{2,\alpha}_m(S^{n-1})$ such that $N(\epsilon, \Lambda_\epsilon, \varphi_\epsilon) = 0$ (i.e. (\ref{phi}) is a positive solution of (\ref{foutball})). Moreover the following estimation holds :
	\[
	|\Lambda_\epsilon| + \left\|\varphi_\epsilon \right\|_{L^{\infty}(S^{n-1})} \leq c \,\epsilon 
	\]
}

\medskip 

The two intermediate results complete the proof of the Proposition.
\hfill $\Box$

\medskip

Observe that, beeing the problem of finding eigenfunctions linear, we can consider as $\phi_\epsilon$ in dimension 2 the following function
\[
\phi_\epsilon = \log \epsilon\, \left[\phi_0 - (\log \epsilon)^{-1}\, (\phi_0(p) + \Lambda_\epsilon)\, \Gamma_p + w_\epsilon + \chi\, H_{\varphi_\epsilon, \epsilon}\right]
\]
This will semplify our reasoning because that function, considered in the coordinates $y = \epsilon x$, converges, in a sense to be made precise, to the function $-\phi_0(p)\, \log |y|$ when $\epsilon$ tends to 0.

\section{Perturbing the complement of a ball}

The following result follows from the implicit function theorem.
\begin{Proposition}\label{perturbing}
\label{pr:1.2}
Given a point $p \in M$, there exists $\epsilon_{0} >0$ and for all $\epsilon \in (0, \epsilon_0)$ and all function $\bar v \in C^{2,\alpha}(S^{n-1})$ satisfying 
\[
\| \bar v \|_{\mathcal C^{2, \alpha} (S^{n-1})}  \leq \epsilon_0 \, ,
\]
and 
\[
\int_{S^{n-1}} \bar v \, \textnormal{ dvol}_{\mathring g} =0 \, ,
\]
there exists a unique positive function $\phi = \phi (\epsilon, p, \bar v) \in {\mathcal C}^{2, \alpha} (B_{1+v}^{g} (p))$, a constant $\lambda = \lambda (\epsilon, p, \bar v) \in \mathbb R$ and a constant $v_0 = v_0 (\epsilon, p, \bar v) \in \mathbb R$ such that
\[
\textnormal{ Vol}_{g} (B_{\epsilon(1+v)}^g(p) ) = \textnormal{ Vol}_{\mathring g} (\mathring B_{\epsilon} )
\]
where $v : =  v_0 + v$ and $\phi$ is a solution to the problem
\begin{equation}\label{formula}
\left\{
\begin{array}{rcccl}
	\Delta_{g} \, \phi + \lambda \, \phi & = & 0 & \textnormal{in} & M \setminus B_{\epsilon(1+v)}^{g}(p) \\[3mm]
	\phi & = & 0 & \textnormal{on} & \partial B_{\epsilon(1+v)}^{g}(p)
\end{array}
\right.
\end{equation}
which is normalized by setting
\begin{equation}
\int_{M \setminus B_{\epsilon(1+v)}^{g} (p)} \phi^2\,  \textnormal{ dvol}_{g}=1 .
\label{noral}
\end{equation} 
In addition $\phi$, $\lambda$ and $v_0$ depend smoothly on the function $\bar v$ and the parameter $\epsilon$. 
\end{Proposition}

\medskip 

{\bf Proof :} We begin by proving that given a point $p \in M$, there exists $\epsilon_{0} >0$ and for all $\epsilon \in (0, \epsilon_0)$ and all function $\bar v \in C^{2,\alpha}(S^{n-1})$ satisfying 
\[
\| \bar v \|_{\mathcal C^{2, \alpha} (S^{n-1})}  \leq \epsilon_0 \, ,
\]
and 
\[
\int_{S^{n-1}} \bar v \, \textnormal{ dvol}_{\mathring g} =0 \, ,
\]
there exists a unique constant $v_0 = v_0 (\epsilon, p, \bar v) \in \mathbb R$ such that
\begin{equation}\label{volume}
\textnormal{ Vol}_{g} (B_{\epsilon(1+v)}^g(p) ) = \textnormal{ Vol}_{\mathring g} (\mathring B_{\epsilon} ) = \epsilon^n\, \textnormal{ Vol}_{\mathring g} (\mathring B_{1} )
\end{equation}
where $v : =  v_0 + \bar v$. Let us define the dilated metric $\bar g = \epsilon^{-2}\, g$. 
Instead of working on a domain depending on the function $v = v_0 + \bar v$, it will be more convenient to work on a fixed domain 
\[
\mathring B_1 :=\{ y \in \mathbb R^n \quad : \quad |y| < 1 \},
\] 
endowed with a metric depending on the function $v$. This can be achieved by considering the  parameterization of $B^{g}_{\epsilon(1+v)}(p) = B^{\bar g}_{(1+v)}(p)$ given by
\[
Y ( y) : =\mbox{Exp}_p^{\bar g} \left( \left(1 + v_0 + \chi ( y) \, \left(\bar v \left(\frac{y}{|y|} \right) \right) \right) \, \sum_i y^i \, E_i \right)
\]
where $\chi$ is a cutoff function identically equal to $0$ when $|y| \leq 1/2$ and identically equal to $1$ when $|y|\geq 3/4$. 

\medskip

Hence (using the result of Proposition~\ref{pr:1.1}) the coordinates we consider from now on are $y \in \mathring B_1$ and in these coordinates the metric $\hat g : = Y^* \bar g$ can be written as 
\[
\hat g  = (1+ v_0)^2 \, \left( \mathring g +  \sum_{i,j} C^{ij} \, dy_i \, dy_j \right) \, ,
\]
where the coefficients $C^{ij} \in {\mathcal C}^{1, \alpha} (\mathring B_\epsilon)$ are functions of $y$ depending on $\epsilon$ $v =v_0+\bar v$ and the first partial derivatives of $v$. Moreover, $C^{ij} \equiv 0$ when $\epsilon =0$ and $\bar v =0$.

\medskip

Observe that 
\[
(\epsilon, v_0, \bar v) \longmapsto C^{ij} (\epsilon, v) \, ,
\]
are smooth maps. 

\medskip

Condition (\ref{volume}), when $\epsilon$ is small enough and not zero, is equivalent to
\begin{equation}
\label{formula-new-2}
\textnormal{ Vol}_{\hat g}(\mathring B_1) =\textnormal{ Vol}_{\mathring g} \, (\mathring B_1)
\end{equation}
that makes sense also for $\epsilon = 0$.
When $\epsilon =0$ and $\bar v \equiv 0$, the metric $\hat g= (1+v_0)^2 \, \mathring g$ is nothing but the Euclidean metric.
We define 
\[
N (\epsilon, \bar v, v_0) : =   \textnormal{ Vol}_{\hat g}(\mathring B_1) - \textnormal{ Vol}_{\mathring g} \, (\mathring B_1)
\]
Observe that $N$ also depends on the choice of the point $p \in M$.

\medskip

We have 
\[
N (0,0,0) ={0}.
\]
It should be clear that the mapping $N$ is a smooth map from a neighborhood of $(0,0,0)$ in  $[0, \infty) \times {\mathcal C}_{m}^{2, \alpha} (S^{n-1}) \times \mathbb R$ into a neighborhood of $0$ in $\mathbb R$.

\medskip

We claim that the partial differential of $N$ with respect to $v_0$, computed at  $(0,0,0,0)$, is given by 
\[
\partial_{v_0} N (0,0,0) = n \, \textnormal{ Vol}_{\mathring g} (\mathring B_1)
\]
Indeed, this time we have $\hat g = (1+ v_0)^{2} \, \mathring g$ since $\bar v \equiv 0$ and $\epsilon =0$ and hence 
\[
N ( 0, 0 , v_0)  =  ((1+ v_0)^{n}-1) \,  \textnormal{ Vol}_{\mathring g}(\mathring B_1) \\[3mm]
\]
So we get 
\[
{\partial_{v_0} N (0,0,0)} =n\,  \textnormal{ Vol}_{\mathring g}(\mathring B_1)  
\]
The claim then follows at once. 

\medskip

Hence the partial differential  of $N$ with respect to both $\psi$ and $v_0$, computed at $(0,0,0)$ is precisely invertible from $\mathbb R$ into $\mathbb R$ and the implicit function theorem ensures, for all {$(\epsilon, \bar v)$} in a neighborhood of $(0,0)$ in $[0, \infty) \times {\mathcal C}^{2, \alpha}_m (S^{n-1})$,  the existence of  a (unique) $v_0 \in \mathcal \mathbb R$ such that $N(\epsilon, \bar v , v_0)=0$. The fact that $v_0$ depends smoothly on the parameter $\epsilon$ and the function $\bar v$ is standard. 

\medskip

Now that we have, for all $0 < \epsilon < \epsilon_0$ and all function $\bar v$ of mean 0, a function $v = v(\epsilon,p,\bar v) \in \mathcal{C}^{2, \alpha}(S^{n-1})$ such that 
\[
\textnormal{ Vol}_{\bar g} (B^{\bar g}_{1+v} ) = \textnormal{ Vol}_{\mathring g} (\mathring B_{1} )
\]
it is easy to find a solution $(\bar \phi, \bar \lambda)$ to the problem (\ref{formula}) and to multiply it by a constant in order to verify the normalization condition. The fact that $\bar \phi$ and $\bar \lambda$ depend smoothly on the parameter $\epsilon$ and the function $\bar v$ is standard.

\hfill $\Box$
\medskip

We will denote the function $\phi = \phi(\epsilon,p,\bar v)$ found in the previous proof as $\phi_{\epsilon,\bar v}$, without noting the dependence on the point $p$. The same for the eigenvalue : $\lambda = \lambda_{\epsilon,\bar v}$, and $\bar \lambda = \bar \lambda_{\epsilon,\bar v} = \epsilon^2\, \lambda_{\epsilon,\bar v}$.
Let us denote
\[
\hat \phi = \hat \phi_{\epsilon,\bar v} = Y^* \phi_{\epsilon,\bar v}
\]
in a neighborhood of $\partial B_\epsilon^g(p)$.
We keep the same notation over all the following paragraphs : for a general $f$ considered in a neighborhood of $\partial B^g_\epsilon(p)$ we will denote
\[
\hat f = Y^* f
\]
We define the operator $F$~:
\[
F (p, \epsilon, \bar v) =  \displaystyle  \hat g (\nabla \hat \phi ,\hat \nu )|_{\partial \mathring B_1}   - \frac{1}{\textnormal{ Vol}_{\hat g} (\partial \mathring B_1)} \, \int_{\partial \mathring B_1} \, \hat g(\nabla \hat \phi , \hat \nu )\, \mbox{dvol}_{\hat g} \, ,
\]
where $\hat \nu$ denotes the unit normal vector field to $\partial \mathring B_1$ and {$(\phi, v_0)$} is the solution of (\ref{formula}) provided by the previous result. Recall that $v = v_0 + \bar v$. Schauder's estimates imply that {$F$} is well defined from a neighborhood of $M \times (0,0)$ in  $M \times [0,\infty) \times \mathcal C^{2,\alpha}_m (S^{n-1})$ into $\mathcal C_m^{1,\alpha}(S^{n-1})$. {Our aim is to find $(p,\epsilon,\bar v)$ such that $F(p,\epsilon,\bar v)=0$. Observe that, with this condition, $\phi$ will be the solution to the problem (\ref{a1}).}

\section{Some estimates}
We want now to give some estimates on $F(p,\epsilon,0)$. In other words we are considering the case when $\bar v = 0$. We keep the notations of the proof of the previous result. If in addition  $v_0 =0$, we can estimate
\[
\hat{g}_{ij}  = \delta_{ij}  + {\mathcal O}(\epsilon^{2}) \, ,
\]
hence 
\[
N(\epsilon, 0, 0)  = {\mathcal O}(\epsilon^{2}) \, .
\]
The implicit function theorem immediately implies that the solution of 
\[
N(\epsilon, 0, v_0) =0
\]
satisfies
\[
|v_0 (\epsilon, p, 0)|\leq c\, \epsilon^2
\]

\medskip

Let us consider the normal coordinates $x$ around $p$. Using the result of Proposition~\ref{pr:1.1} is possible to show that
\[
\begin{array}{rllll}
g^{ij} & = & \displaystyle \delta_{ij} - \frac{1}{3} \, R_{ikj\ell} \, x^{k} \, x^{\ell} - \frac{1}{6} \, R_{ikj\ell,m} \, x^{k} \, x^{\ell}  \, x^{m} + {\mathcal O} (|x|^{4})\\[3mm]
\log |g| & = & \displaystyle \frac{1}{3} \, R_{k\ell} \, x^{k} \, x^{\ell}  + \frac{1}{6} \, R_{k\ell,m} \, x^{k} \, x^{\ell} \, x^{m} + {\mathcal O} (|x|^{4})
\end{array}
\]where
$$
R_{k\ell} = \sum_{i=1}^{n} R_{iki\ell}\ \textnormal{and}\ R_{k\ell,m} = \sum_{i=1}^{n} R_{iki\ell,m}
$$

\medskip

A straightforward calculation allows us to obtain the expansion of $\Gamma_p$.
Recall that 
\[
\Delta_{g}  : =  \sum_{i,j} g^{ij} \, \partial_{x_i} \partial_{x_j} + \sum_{i,j} \partial_{x_i}  g^{ij} \, \partial_{x_j} + \frac{1}{2} \, \sum_{i,j} g^{ij} \, \partial_{x_i} \log |g| \, \partial_{x_j}  
\]
From the definition of the Green function $\Gamma_p$ we can obtain its expansion near $p$. For $n \geq 5$ we have
\begin{equation}\label{gamma5}
\begin{array}{rllll}
 \displaystyle \Gamma_p (x)  & = & \displaystyle |x|^{2-n}\, +\\[5mm]
 & + & \displaystyle \Bigg ( \frac{2-n}{18}\, R_{ikj\ell}\, x^i\, x^k\, x^j\, x^{\ell}\, |x|^{-n}\, -\, \frac{1}{12}R_{j \ell}\, x^j\, x^{\ell}\, |x|^{2-n}\, + \frac{\textnormal{Scal}(p) - 6\,\lambda_0}{12(4-n)}\, |x|^{4-n}\Bigg)\,+ \\[5mm] 
 & + & \displaystyle \Bigg( \frac{2-n}{48}\, R_{ikj\ell, t}\, x^i\, x^k\, x^j\, x^{\ell}\, x^t\, |x|^{-n}\, + 
 \frac{1}{36}\,R_{\cdot kj\ell, \cdot}\, x^k\, x^j\, x^{\ell}\, |x|^{2-n}\, +\\[5mm]
 & & \displaystyle \qquad -\, \frac{1}{24}\, R_{j\ell, t}\, x^j\, x^{\ell}\, x^t\, |x|^{2-n} + \frac{3\,{\textnormal{Scal}}_{,t}}{64(4-n)}\, x^t\, |x|^{4-n}\Bigg)\, + a + \\[5mm]
 & + &  \mathcal{O}(|x|^{6-n}).
\end{array}
\end{equation}
When $n = 4$ we have
\begin{equation}\label{gamma4}
\begin{array}{rllll}
 \displaystyle \Gamma_p (x)  & = & \displaystyle |x|^{-2}\, +\\[7mm]
 & + & \displaystyle \Bigg ( -\frac{1}{9}\, R_{ikj\ell}\, x^i x^k x^j x^{\ell}\, |x|^{-4}\, -\, \frac{1}{12}R_{j \ell}\, x^j x^{\ell}\, |x|^{-2}\, + \frac{\textnormal{Scal}(p) - 6\,\lambda_0}{12}\, \log|x|\Bigg)\,+ \\[7mm] 
 & + & \displaystyle \Bigg( -\frac{1}{24}\, R_{ikj\ell, t}\, x^i x^k x^j x^{\ell} x^t\, |x|^{-4}\, + 
 \frac{1}{36}\,R_{\cdot kj\ell, \cdot}\, x^k x^j x^{\ell}\, |x|^{-2}\, +\\[7mm]
 & & \displaystyle \qquad -\, \frac{1}{24}\, R_{j\ell, t}\, x^j x^{\ell} x^t\, |x|^{-2} + \frac{3\,{\textnormal{Scal}}_{,t}}{64}\, x^t\, \log |x|\Bigg)\, +\\[7mm]
 & + & a' + b \cdot x + \mathcal{O}(|x|^{\alpha}),
\end{array}
\end{equation}
for all $\alpha < 2$, where $a, a'$ are  constants and $b \in \mathbb{R}^n$.
In the above expressions we used the notation
$$
R_{\cdot kj\ell,\cdot} : = \sum_{i=1}^{n} R_{ikj\ell,i}
$$
Observe that to find such an expression we used the fact that $R(X,X)\equiv 0$, the symmetries of the curvature tensor for which if either $i=k$ or $j=\ell$ then $R_{ikj\ell,m} = 0$, and the second Bianchi identity 
\[
\sum_{j} R_{t j,j} = \sum_{j} R_{j t,j} = \frac{1}{2} \, \textnormal{ Scal}_{,t}.
\]
Remark that for $n=2$ or $n=3$ is not possible to give a precise expansion of the Green function near $p$ using only the local part of the equation that defines $\Gamma_p$. 

\medskip 

The main result of this section is the~:
\begin{Proposition}
\label{le:3.3}
\medskip
In the CASE 1 (i.e. the case where $\phi_0$ is not constant) there exists a constant $c >0$ such that, for all $p \in M$ and all  $\epsilon \geq 0$ small enough we have 
\[
\begin{array}{lll}
\| F (p, \epsilon, 0) \|_{\mathcal C^{1, \alpha}} \leq c \, \epsilon &  \textnormal{if} & n \geq 3\\[3mm]
\| F (p, \epsilon, 0) \|_{\mathcal C^{1, \alpha}} \leq c \, \epsilon\, \log \epsilon &  \textnormal{if} & n =2
\end{array}
\] 
Moreover there exists a constant $C_n$ depending only on $n$, such that for all $a \in \mathbb R^n$ the following estimates hold
\[
\begin{array}{lll}
\displaystyle \left| \int_{S^{n-1}} \mathring g (a, \cdot ) \, F(p,\epsilon,0) \, \textnormal{ dvol}_{\mathring g} - C_n \,\epsilon\, g( \nabla \phi_0 (p) , \Theta (a) ) \right| \leq c \,  \epsilon^2 \, \| a\| &  \textnormal{if} & n \geq 3\\[7mm]
\displaystyle \left| \int_{S^{n-1}} \mathring g (a, \cdot ) \, F(p,\epsilon,0) \, \textnormal{ dvol}_{\mathring g} - C_n \,\epsilon\,  \log \epsilon\, g( \nabla \phi_0 (p) , \Theta (a) ) \right| \leq c \,  \epsilon^2\, \log \epsilon \, \| a\| &  \textnormal{if} & n = 2
\end{array} 
\]
In the CASE 2 (i.e. the case where $\phi_0$ is a constant function) and for $n \geq 4$ there exists a constant $c >0$ such that, for all $p \in M$ and all  $\epsilon \geq 0$ small enough we have 
\[
\begin{array}{lll}
\| F (p, \epsilon, 0) \|_{\mathcal C^{1, \alpha}} \leq c \, \epsilon^2 &  \textnormal{if} & n \geq 5\\[3mm]
\| F (p, \epsilon, 0) \|_{\mathcal C^{1, \alpha}} \leq c \, \epsilon^2\, \log \epsilon &  \textnormal{if} & n =4
\end{array}
\] 
Moreover there exists a constant $C_n$ (depending only on $n$), such that for all $a \in \mathbb R^n$ the following estimates hold:
\[
\begin{array}{rll}
\displaystyle \left| \int_{S^{n-1}} \mathring g (a, \cdot ) \, F(p,\epsilon,0) \, \textnormal{ dvol}_{\mathring g} - C_n \, \epsilon^{3} \, g( \nabla \textnormal{ Scal} (p) , \Theta (a) ) \right| \leq c \,  \epsilon^{4} \, \| a\| &  \textnormal{if} & n \geq 5\\[7mm]
\displaystyle \left| \int_{S^{n-1}} \mathring g (a, \cdot ) \, F(p,\epsilon,0) \, \textnormal{ dvol}_{\mathring g} - C_n \, \epsilon^{3} \log \epsilon \, g( \nabla \textnormal{ Scal} (p) , \Theta (a) ) \right| \leq c \,  \epsilon^{3} \, \| a\| &  \textnormal{if} & n = 4
\end{array} 
\]
\end{Proposition}

{\bf Proof :} Let $\epsilon$ be small enough, and $\bar v = 0$. We know that $v_0 = \mathcal{O}(\epsilon^2)$, then from proposition \ref{eigenfunction} it follows that for all $\epsilon$ small enough there exists $(\Lambda_\epsilon,\varphi_\epsilon,w_\epsilon)$ in a neighborhood of $(0,0,0)$ in $\mathbb{R} \times \mathcal{C}^{2,\alpha}_m(S^{n-1}) \times \mathcal{C}^{2,\alpha}(M \setminus B_\epsilon^g(p))$ such that the first eigenfunction of $-\Delta_g$ over the complement of $B_{\epsilon(1+v_0)}^g(p)$ with 0 Dirichlet condition at $\partial B_{\epsilon(1+v_0)}^g(p)$ is given by
\[
\phi_{\epsilon,0} = \phi_0 - \epsilon^{n-2}\, (1+v_0)^{n-2}\,(\phi_0(p) + \Lambda_\epsilon)\, \Gamma_p + w_\epsilon + \chi\, H_{\varphi_\epsilon, \epsilon}
\]
if $n \geq 3$, and by
\[
\phi_{\epsilon,0} = \log (\epsilon\,(1+v_0))\, \left[\phi_0 - (\log \epsilon\,(1+v_0))^{-1}\, (\phi_0(p) + \Lambda_\epsilon)\, \Gamma_p + w_\epsilon + \chi\, H_{\varphi_\epsilon, \epsilon}\right]
\]
if $n = 2$, where estimations given in propositions \ref{eigenfunction} and \ref{eigenfunction2} continue to hold because $v_0 = \mathcal{O}(\epsilon^2)$. 

\medskip 
 
From the expression of $\phi_{\epsilon,0}$ it follows that in the CASE 1 we have
\begin{eqnarray*}
\int_{S^{n-1}} \mathring g (a, \cdot ) \,\hat g (\nabla \hat \phi_{\epsilon,0},\hat \nu)|_{\partial \mathring B_1}\,\textnormal{ dvol}_{\mathring g}&  = & (1 + \mathcal{O}(\epsilon)) \, \int_{S^{n-1}} \mathring g (a, \cdot ) \, \frac{\partial \hat \phi_{\epsilon,0}}{\partial r}|_{\partial \mathring B_1} \, \textnormal{ dvol}_{\mathring g} \\[3mm]
& = & (1 + \mathcal{O}(\epsilon)) \, \int_{S^{n-1}} \mathring g (a, \cdot ) \, \frac{\partial \hat \phi_0}{\partial r}|_{\partial \mathring B_1} \, \textnormal{ dvol}_{\mathring g} + \mathcal{O}(\epsilon^2)\\[3mm]
& = & C_n \, \epsilon \, g( \nabla \phi_0 (p) , \Theta (a) ) + \mathcal{O}(\epsilon^2)
\end{eqnarray*}
for $n \geq 3$, and
\[
\int_{S^{n-1}} \mathring g (a, \cdot ) \,\hat g (\nabla \hat \phi_{\epsilon,0},\hat \nu)|_{\partial \mathring B_1}\,\textnormal{ dvol}_{\mathring g}  =  C_2 \,\epsilon\, \log \epsilon\, g( \nabla \phi_0 (p) , \Theta (a) ) + \mathcal{O}(\epsilon^2\, \log \epsilon)
\]
for $n = 2$.
Then all the estimates for the CASE 1 follow at once from this computation together with the fact that, when $\bar v \equiv 0$, the unit normal vector $\hat \nu$ about the boundary is given by $(1 + v_0)\, \partial_r\, (1 + \mathcal{O}(\epsilon))$ because the metric $\hat g$ near $p$ is the euclidean metric multiplied by $(1+v_0)^2$ and perturbed by some $\mathcal{O}(\epsilon^2)$ terms.

\medskip 

For the CASE 2 the situation is much more complex. We remark that if $\phi_0$ is constant, then 
\[
\int_{S^{n-1}} \mathring g (a, \cdot ) \, \hat g (\nabla \hat \phi_0, \hat \nu)|_{\partial \mathring B_1}\,\textnormal{ dvol}_{\mathring g} = 0.
\]

\medskip

Let us compute now 
\[
\epsilon^{n-2}\,(1+v_0)^{n-2}\,(\phi_0(p) + \Lambda_\epsilon)\, \int_{S^{n-1}} \mathring g (a, \cdot ) \, \hat g (\nabla \hat \Gamma_p, \hat \nu)|_{\partial \mathring B_1}\,\textnormal{ dvol}_{\mathring g}
\]
We remark that the previous term is equal to 
\[
(1 + \mathcal{O}(\epsilon))\, \epsilon^{n-2}\,\phi_0(p)\, \int_{S^{n-1}} \mathring g (\cdot,  a)  \frac{\partial \hat \Gamma_p}{\partial r} \, \textnormal{ dvol}_{\mathring g}
\]
For this reason we will compute 
\[
\epsilon^{n-2}\, \phi_0(p)\, \int_{S^{n-1}} \mathring g (\cdot,  a)  \frac{\partial \hat \Gamma_p}{\partial r} \, \textnormal{ dvol}_{\mathring g} 
\] 
Recall that
\[
\hat \Gamma_p(y) = \Gamma_p (\epsilon\,(1+v_0)\, y)
\]
in a neighborhood of $\partial \mathring B_1$, then from (\ref{gamma5}) and (\ref{gamma4}) (keeping in mind that $v_0 = \mathcal{O}(\epsilon^2)$) we obtain easily the expression of $\hat \Gamma_p(y)$ in power of $\epsilon$.
Observe that, in the expansion of $\hat \Gamma_p$, terms which contain an even number of coordinates, such as $y^{i}y^{j}y^{k}y^{\ell}$ or $y^{j}y^{\ell}$ etc. do not contribute to the result since, once derived with respect to $r$ they continue to contain an even number of coordinates, and multiplied then by $\mathring g(y , a)$, their average over $S^{n-1}$ is $0$. Then, considering only terms which contain an odd number of coordinates we have for $n \geq 5$:
\[
\begin{array}{l}
 \displaystyle \epsilon^{n-2}\, \int_{S^{n-1}} \mathring g  (y, a)\, \frac{\partial \hat \Gamma_p}{\partial r}\, \mbox{dvol}_{\mathring g}  =\\[7mm]
   \qquad \displaystyle = \epsilon^{3}\,a_{\sigma} \Bigg[\int_{S^{n-1}} y^{\sigma}\, \cdot \frac{y^{\tau}}{|y|}\,\cdot \frac{\partial}{\partial y^{\tau}} \Bigg(\frac{2-n}{48}\, R_{ikj\ell, t}\, y^i y^k y^j y^{\ell} y^t\, |y|^{-n} + \frac{1}{36}\,R_{\cdot kj\ell, \cdot}\, y^k y^j y^{\ell}\, |y|^{2-n}\\[7mm]
   \qquad \qquad + \displaystyle \frac{3 {\textnormal{Scal}}_{,t}}{64(4-n)}\, y^t\, |y|^{4-n} -\, \frac{1}{24} \, R_{j\ell, t}\, y^j y^{\ell} y^t\, |y|^{2-n} \Bigg) \mbox{dvol}_{\mathring g} \Bigg] + {\mathcal O} (\epsilon^{4})\\[7mm]
  \qquad = \displaystyle \epsilon^{3}\, (5-n)\, a_{\sigma} \Bigg[\int_{S^{n-1}} y^{\sigma}\, \Bigg(\frac{2-n}{48}\, R_{ikj\ell, t}\, y^i y^k y^j y^{\ell} y^t + \frac{1}{36}\,R_{\cdot kj\ell, \cdot}\, y^k y^j y^{\ell} + \frac{3{\textnormal{Scal}}_{,t}}{64(4-n)}\, y^t\\ [7mm]
 \displaystyle \qquad \qquad -\, \frac{1}{24} \, R_{j\ell, t}\, y^j y^{\ell} y^t \Bigg) \mbox{dvol}_{\mathring g} \Bigg] + {\mathcal O} (\epsilon^{4})
\end{array}
\] 

We make use of the identities in the Appendix to conclude that there exists a costant $C_n^{(1)}$ such that 
\begin{equation}\label{abc}
 \displaystyle \epsilon^{n-2}\, \phi_0(p)\, \int_{S^{n-1}} \mathring g  (y, a) \, \frac{\partial \hat \Gamma_p}{\partial r}(y) =  C_n^{(1)} \epsilon^{3} \, {g \big(\nabla \textnormal{ Scal(p)},  \Theta (a) \big)}\, +\, \mathcal{O}(\epsilon^{4}).
\end{equation}
where we have
\[
C_n^{(1)} = \frac{5-n}{4n}\,{\textnormal{Vol}}_{\mathring g}(S^{n-1})\, \left[ -\frac{1}{3(n+2)}+\frac{3}{16(4-n)}\right]\, \phi_0(p)
\]
Remark that when $n= 5$ such a constant is 0. For $n = 4$ we have
\[
\begin{array}{l}
 \displaystyle \epsilon^{n-2}\, \int_{S^{n-1}} \mathring g  (y, a)\, \frac{\partial \hat \Gamma_p}{\partial r}\, \mbox{dvol}_{\mathring g}  =\\[7mm]
   \qquad \displaystyle = \epsilon^{3}\,\log \epsilon\, a_{\sigma} \Bigg[\int_{S^{3}} y^{\sigma}\, \cdot \frac{y^{\tau}}{|y|}\,\cdot \frac{\partial}{\partial y^{\tau}} \left(\frac{3 {\textnormal{Scal}}_{,t}}{64}\, y^t\, \log|y|\right)\,  \mbox{dvol}_{\mathring g} \Bigg] + {\mathcal O} (\epsilon^{3})\\[7mm]
  \qquad = \displaystyle \frac{3}{256}\, \textnormal{Vol}_{\mathring g}(S^{3})\, \epsilon^{3}\, \log \epsilon\, g \big(\nabla \textnormal{ Scal(p)},  \Theta (a) \big)\, + {\mathcal O} (\epsilon^{3})
\end{array}
\] 
and then we set $C^{(1)}_4 = \displaystyle \frac{3}{256}\, \textnormal{Vol}_{\mathring g}(S^{3})$.

 \medskip
 
The last term we have to compute is
\[
\int_{S^{n-1}} \mathring g (a, \cdot ) \, \hat g (\nabla (\hat w_\epsilon + \hat H_{\varphi_\epsilon,\epsilon}), \hat \nu)|_{\partial \mathring B_1}\,\textnormal{ dvol}_{\mathring g}
\]
As before we have
\[
\begin{array}{l}
\displaystyle \int_{S^{n-1}} \mathring g (a, \cdot ) \, \hat g (\nabla (\hat w_\epsilon + \hat H_{\varphi_\epsilon,\epsilon}), \hat \nu)|_{\partial \mathring B_1}\,\textnormal{ dvol}_{\mathring g} =\\[3mm]
\qquad \qquad \displaystyle = (1 + \mathcal{O}(\epsilon))\, \int_{S^{n-1}} \mathring g (a, \cdot ) \, \frac{\partial (\hat w_\epsilon + \hat H_{\varphi_\epsilon,\epsilon})}{\partial r}|_{\partial \mathring B_1}\,\textnormal{ dvol}_{\mathring g} 
\end{array}
\]
In the Proposition \ref{eigenfunction} we proved that in the CASE 2 
\[
\|w_\epsilon\|_{\mathcal{C}^{2,\alpha}_\nu(M \setminus \{p\})} \leq c\, \left(\epsilon^{4} + \epsilon^{\beta-\nu}\right)\,\,\, \forall \beta < 4
\]
for $n=4$ and 
\[
\|w_\epsilon\|_{\mathcal{C}^{2,\alpha}_\nu(M \setminus \{p\})} \leq c\, \left(\epsilon^{2n-4} + \epsilon^{1+n} + \epsilon^{4-\nu}\right)
\]
for $n \geq 5$.
Hence
\[
\|\nabla \hat w_\epsilon\|_{L^{\infty}(\partial \mathring B_1)} \leq c\, (\epsilon^{4+\nu} + \epsilon^{\beta})\,\,\, \forall \beta < 4
\]
for $n=4$ and 
\[
\|\nabla \hat w_\epsilon\|_{L^{\infty}(\partial \mathring B_1)} \leq c\, (\epsilon^{2n-4+\nu} + \epsilon^{1+n+\nu} + \epsilon^{4})
\]
for $n \geq 5$
(keep in mind that we are estimating the gradient of the dilated function $\hat w_\epsilon$). Remember that $\nu \in (2-n,4-n)$ because $n \geq 4$. It follows that we can choose $\nu$ in order to have
\[
\int_{S^{n-1}} \mathring g (a, \cdot ) \, \frac{\partial \hat w_\epsilon}{\partial r}|_{\partial \mathring B_1}\,\textnormal{ dvol}_{\mathring g} = \mathcal{O}(\epsilon^\beta)
\]
with $\beta = 4$ for $n \geq 5$ and $\beta = 3$ for $n = 4$. Let us consider now $\hat H_{\varphi_\epsilon,\epsilon}$. We do not know the expression of $\hat H_{\varphi_\epsilon,\epsilon}$ in a neighborhood of $\partial \mathring B_1$, but we can know its value on $\partial \mathring B_1$. From the equality $\hat \phi_\epsilon = 0$ on $\partial \mathring B_1$, using the the estimate on the function $\hat w_{\epsilon}$, we have that 
\[
\begin{array}{rllll}
 \displaystyle \hat H_{\varphi,\epsilon}  & = & \displaystyle - \phi_0(p) + (1+v_0)^{n-2}\, (\phi_0(p) + \Lambda_\epsilon)\, \Bigg[\epsilon^2\, \Bigg ( \frac{2-n}{18}\, R_{ikj\ell}\, y^i y^k y^j y^{\ell}\, -\, \frac{1}{12}R_{j \ell}\, y^j y^{\ell}\, + \\[7mm]
 & & \displaystyle \qquad + \frac{\textnormal{Scal}(p) - 6 \lambda_0}{12\, (4-n)}\Bigg) + \epsilon^3\, \Bigg( \frac{2-n}{48}\, R_{ikj\ell, t}\, y^i y^k y^j y^{\ell} y^t\,  + 
 \frac{1}{36}\,R_{\cdot kj\ell, \cdot}\, y^k y^j y^{\ell}\,+\\[7mm]
 & & \qquad \qquad \displaystyle -\, \frac{1}{24}\, R_{j\ell, t}\, y^j y^{\ell} y^t\, + \frac{3\,{\textnormal{Scal}}_{,t}}{64(4-n)}\, y^t\, \Bigg)\Bigg] + \mathcal{O}(\epsilon^4).
\end{array}
\]
on $\partial \mathring B_1$, for $n \geq 5$. 
For $ n = 4$ :
\[
\begin{array}{rllll}
 \displaystyle \hat H_{\varphi,\epsilon}  & = & \displaystyle -\phi_0(p) + (1+v_0)^{n-2}\, (\phi_0(p) + \Lambda_\epsilon)\, \Bigg[\epsilon^2\,\log \epsilon\, \frac{\textnormal{Scal}(p) - 6 \lambda_0}{12}\, + \\[7mm]
 & & \displaystyle \qquad \epsilon^2\, \Bigg ( -\frac{1}{9}\, R_{ikj\ell}\, y^i y^k y^j y^{\ell}\, -\, \frac{1}{12}R_{j \ell}\, y^j y^{\ell}\Bigg)\, +\, \epsilon^3\, \log \epsilon\,  \frac{3\,{\textnormal{Scal}}_{,t}}{64}\, y^t\, \Bigg)\Bigg] + \mathcal{O}(\epsilon^3).
\end{array}
\]
Let us define an harmonic extension of $\mathring g (y,a)$ to $\mathbb{R}^{n} \setminus \mathring B_1$ :
\begin{equation}\label{Hharmonic}
\left\{
\begin{array}{rclll}
	\Delta_{\mathring g} G_{a} & = & 0 & \textnormal{in} & \mathbb{R}^{n} \setminus \mathring B_1\\[3mm]
	 G_{a} & = & \mathring g (y,a) & \textnormal{on} & \partial \mathring B_1 
\end{array}
\right.
\end{equation}
It is easy to check that
\[
G_a(y) = |y|^{-n}\, \mathring g (y,a)
\]
We observe that the functions $G_a$ and $\hat H_{\varphi_\epsilon,\epsilon}$, by Lemma \ref{lemmaharmonic} converge to $0$ when $|y| \rightarrow + \infty$. Then
\[
\begin{array}{lll}
\displaystyle \int_{S^{n-1}} \mathring g (a, \cdot ) \, \frac{\partial \hat H_{\varphi_\epsilon,\epsilon}}{\partial r}|_{\partial \mathring B_1}\,\textnormal{ dvol}_{\mathring g} & = & \displaystyle \int_{S^{n-1}} \hat H_{\varphi_\epsilon,\epsilon}\, \frac{\partial G_a}{\partial r}|_{\partial \mathring B_1}\,\textnormal{ dvol}_{\mathring g} = \\[3mm]
 & & \qquad \displaystyle = (1-n) \int_{S^{n-1}} \hat H_{\varphi_\epsilon,\epsilon}\, \mathring g (y,a)\,\textnormal{ dvol}_{\mathring g}
\end{array}
\]
Using the expansion of $\hat H_{\varphi_\epsilon,\epsilon}$ that we found with the identities in the Appendix, we conclude that there exists a constant $C_n^{(2)} \neq 0$ such that 
\begin{equation}
 \displaystyle \int_{S^{n-1}} \mathring g  (y, a) \,\hat H_{\varphi_\epsilon,\epsilon}\,\textnormal{ dvol}_{\mathring g}  =  C_n^{(2)}\, \epsilon^{3} \, {g \big(\nabla \textnormal{ Scal(p)},  \Theta (a) \big)}\, +\, \mathcal{O}(\epsilon^{4}).
\end{equation}
where 
\[
C_n^{(2)} = \frac{1}{4n}\,{\textnormal{Vol}}_{\mathring g}(S^{n-1})\, \left[ -\frac{1}{3(n+2)}+\frac{3}{16(4-n)}\right]\, \phi_0(p)
\]
for $n \geq 5$, and for $n=4$
\begin{equation}
 \displaystyle \int_{S^{n-1}} \mathring g  (y, a) \,\hat H_{\varphi_\epsilon,\epsilon}\,\textnormal{ dvol}_{\mathring g}  =  C_4^{(2)}\, \epsilon^{3} \log \epsilon \, {g \big(\nabla \textnormal{ Scal(p)},  \Theta (a) \big)}\, +\, \mathcal{O}(\epsilon^{3}).
\end{equation}
with
\[
C_4^{(2)} = \frac{3}{256}\, \textnormal{Vol}_{\mathring g}(S^{3}).
\] 

\medskip 

Summarizing we conclude that in the CASE 2
\[
\| F (p, \epsilon, 0) \|_{\mathcal C^{1, \alpha}} = \mathcal{O}(\epsilon^2)
\] 
and there exists a constant $C_n$ depending only on $n$, such that for all $a \in \mathbb R^n$ the following estimates hold : for $n \geq 5$
\[
\left| \int_{S^{n-1}} \mathring g (a, \cdot ) \, F(p,\epsilon,0) \, \textnormal{ dvol}_{\mathring g} - C_n \, \epsilon^{3} \, g( \nabla \textnormal{ Scal} (p) , \Theta (a) ) \right| \leq c \,  \epsilon^{4} \, \| a\| \, .
\]
where
\[
C_n = \frac{6-2n}{n}\,{\textnormal{Vol}}_{\mathring g}(S^{n-1})\, \left[ -\frac{1}{3(n+2)}+\frac{3}{16(4-n)}\right]\, \phi_0(p)
\]
and for $n = 4$
\[
\left| \int_{S^{n-1}} \mathring g (a, \cdot ) \, F(p,\epsilon,0) \, \textnormal{ dvol}_{\mathring g} - C_4 \, \epsilon^{3} \log \epsilon\, g( \nabla \textnormal{ Scal} (p) , \Theta (a) ) \right| \leq c \,  \epsilon^{3} \, \| a\| \, .
\]
where
\[
C_4 = -\frac{3}{128}\,{\textnormal{Vol}}_{\mathring g}(S^{3})\, \phi_0(p)
\]
Remark that $C_n\neq0$ for all $n \geq 4$. 
This completes the proof of the result.  \hfill $\Box$

\section{Linearizing the operator $F$}

Our next task will be to understand the structure of $L_0$, the operator obtained by linearizing $F$ with respect to $\bar v$ at $\epsilon =0$ and $\bar v=0$. We will see that this operator is a first order elliptic operator which does not depend on the point $p$. 

\medskip

Let us define in $\mathbb{R}^n \setminus \{0\}$
\[
\phi_1(y) = \left\{
\begin{array}{lll}
	\displaystyle \phi_0(p)\, (1-|y|^{2-n}) & \textnormal{if} \qquad n \geq 3\\[3mm]
	\phi_0(p)\, \log|y| & \textnormal{if} \qquad n=2
\end{array}
\right.
\] 
For all $ \bar v \in \mathcal C^{2, \alpha}_m (S^{n-1})$ let $\psi$ be the (unique) bounded solution of 
\begin{equation}\label{c00}
\left\{
\begin{array}{rcll}
	\displaystyle \Delta_{\mathring g}\, \psi & = & 0  & \textnormal{in} \qquad \mathbb R^n \setminus \mathring B_1 \\[3mm]
	\psi & = &  - \displaystyle  {\partial_r \phi_{1}} \, \bar v  & \textnormal{on}\qquad \partial \mathring B_1
\end{array}
\right.
\end{equation}
By the Lemma \ref{lemmaharmonic}, $|\psi(y)| \longrightarrow 0$ when $|y| \rightarrow \infty$. We define
\begin{equation}
H(\bar v) : = \left( {\partial_r \psi }  + {\partial^2_r \phi_{1}} \, \bar v \right) \, |_{\partial \mathring B_1}
\label{Hache}
\end{equation}

\medskip

We will need the following result~: 
\begin{Proposition}
\label{H}
The operator 
\[
H : \mathcal C^{2, \alpha}_m (S^{n-1}) \longrightarrow \mathcal C^{1, \alpha}_m (S^{n-1}) , 
\]
is a self adjoint, first order elliptic operator. The kernel of $H$ is given by $V_1$, the eigenspace of $-\Delta_{S^{n-1}}$ associated to the eigenvalue $n-1$. Moreover there exists $c >0$ such that 
\[
\| w \|_{\mathcal C^{2, \alpha}(S^{n-1})}  \leq c \, \| H(w)  \|_{\mathcal C^{1, \alpha}(S^{n-1})} \, ,
\]
provided $w$ is $L^2(S^{n-1})$-orthogonal to $V_0 \oplus V_1$, where $V_0$ is the eigenspace associated to constant functions.
\end{Proposition}

{\bf Proof :}  The fact that $H$ is a first order elliptic operator is standard since it is the sum of the Dirichlet-to-Neumann operator for $\Delta_{\mathring g}$ and  a constant times the identity. In particular, elliptic estimates yield
\[
\| H(w)  \|_{\mathcal C^{1, \alpha}(S^{n-1})} \leq c \, \| w \|_{\mathcal C^{2, \alpha}(S^{n-1})}
\]

The fact that the operator $H$ is (formally) self-adjoint is easy. Let $\psi_1$ (resp. $\psi_2$) the solution of (\ref{c00}) corresponding to the function ${w_{1}}$ (resp. $w_2$). We compute
\[
\begin{array}{rlllll}
\partial_r \phi_1 (1) \, \displaystyle \int_{\partial \mathring B_1} ( H(w_1) \, w_2 - w_1 \, H(w_2) )    \, \mbox{dvol}_{\mathring g} & = & \partial_r \phi_1 (1) \, \displaystyle \int_{\partial \mathring B_1}( \partial_r \psi_1 \, w_2 - \partial_r \psi_2 \, w_1)  \, \mbox{dvol}_{\mathring g}\\[3mm]
& = & \displaystyle \int_{\partial \mathring B_1}( \psi_1 \, \partial_r \psi_2 - \psi_2  \, \partial_r  \psi_1  )  \, \mbox{dvol}_{\mathring g}\\[3mm]
& = & \displaystyle \lim_{R \rightarrow \infty} \int_{\mathring B_R \setminus \mathring B_1}( \psi_1  \,  \Delta_{\mathring g} \psi_2 - \psi_2 \Delta_{\mathring g} \, \psi_1 )  \, \mbox{dvol}_{\mathring g}\\[3mm]
&  & - \displaystyle \lim_{R \rightarrow \infty} \int_{\partial \mathring B_R}( \psi_1 \, \partial_r \psi_2 - \psi_2  \, \partial_r  \psi_1  )  \, \mbox{dvol}_{\mathring g}\\[3mm]
& = & 0
\end{array}
\]
 
Let us consider
\[
w = \sum_{j \geq 1} w_j
\]
the eigenfunction decomposition of $w$. Namely $w_j \in V_j$. We can define $\psi_j$ to be the bounded solution of 
\begin{equation}
\left\{
\begin{array}{rcll}
	\displaystyle \Delta_{\mathring g}\, \psi_j & = & 0  & \textnormal{in} \qquad \mathbb R^n \setminus \mathring B_1 \\[3mm]
	\psi_j & = &  - \displaystyle  {\partial_r \phi_{1}} \, w_j  & \textnormal{on}\qquad \partial \mathring B_1
\end{array}
\right.
\end{equation}
i.e. 
\[
\psi_j (y)= -|y|^{2-n-j}\, w_j(y/|y|)\, \partial_r \phi_1|_{\partial \mathring B_1}
\]
Then 
\[
H(w) = \sum_j \partial_r \psi_j + \partial_r^2 \phi_1|_{\partial \mathring B_1} \, w = \sum_j \left[-(2-n-j)\, \partial_r \phi_1|_{\partial \mathring B_1} + \partial_r^2 \phi_1|_{\partial \mathring B_1} \right]\, w_j
\]
With this alternative formula for $H$, it is clear that $H$ preserves the eigenspaces $V_j$ and in particular, $H$ maps into the space of functions whose mean over $S^{n-1}$ is 0. Moreover, it is easy to see that $V_1$ is the only kernel of the operator. In fact, 
\[
\partial_r \phi_1|_{\partial \mathring B_1}  = \left\{
\begin{array}{lll}
	\displaystyle -(2-n)\, \phi_0(p) & \textnormal{if} \qquad n \geq 3\\[3mm]
	\phi_0(p) & \textnormal{if} \qquad n=2
\end{array}
\right.
\] 
and
\[
\partial_r^2 \phi_1|_{\partial \mathring B_1}  = \left\{
\begin{array}{lll}
	\displaystyle -(2-n)(1-n)\, \phi_0(p) & \textnormal{if} \qquad n \geq 3\\[3mm]
	-\phi_0(p) & \textnormal{if} \qquad n=2
\end{array}
\right.
\] 
and then $H(w_j) = 0$ if and only if $j = 1$. This completes the proof of the result.
 \hfill  $\Box$

\medskip

The main result of this section is the following~: 
\begin{Proposition}\label{l0}
The operator $L_0$ is equal to $H$.
\end{Proposition}

{\bf Proof~:} By definition, the operator $L_0$ is the linear operator obtained by linearizing $F$ with respect to $\bar v$ at $\epsilon =0$ and $\bar v=0$. In other words, we  have
\[
L_0 (\bar w) = \lim_{s \rightarrow 0} \frac{F(p,0, s\, \bar w) - F(p,0,0)}{s} .
\]
It is easy to see that $F(p,0,0) = 0$. In fact, we saw that the first eigenfunction $ \phi_{\epsilon,0}$ over $M \setminus B_{\epsilon(1+v_0)}^{g}(p)$ is given by
\[
\begin{array}{lll}
\phi_{\epsilon,0} = \phi_0 - \epsilon^{n-2}\,(1 + v_0)^{n-2}\, (\phi_0(p) + \Lambda_\epsilon)\, \Gamma_p + w_\epsilon + \chi\, H_{\varphi_\epsilon, \epsilon} & \textnormal{if} & n \geq 3\\[3mm]
\phi_{\epsilon,0} = \log(\epsilon\,(1+v_0))\, \left[\phi_0 - (\log (\epsilon\,(1+v_0)))^{-1}\,(\phi_0(p) + \Lambda_\epsilon)\, \Gamma_p + w_\epsilon + \chi\, H_{\varphi_\epsilon, \epsilon}\right] & \textnormal{if} & n =2
\end{array}
\] 
where $v_0 = v_0(p,\epsilon,0) = \mathcal{O}(\epsilon^2)$, for some $(\Lambda_\epsilon,w_\epsilon,\varphi_\epsilon) \in \mathbb{R} \times C^{2,\alpha}(M \setminus B^{\bar g}_{1 + v_0}(p)) \times C^{2,\alpha}_m(S^{n-1})$, where the estimations of Proposition \ref{le:3.3} hold because $v_0 = \mathcal{O}(\epsilon^2)$.
If we consider this expressions only in a neighborhood of $\partial B_{\epsilon(1+v_0)}^{g}(p)$ and the parameterization $Y$ given in the proof of the Proposition \ref{perturbing} with coordinates $y$ in a neighborhood of $\partial \mathring B_1$, it is easy to see that the function $\hat \phi_0 = Y^* \phi_0$ is equal to the constant function $\hat \phi_0 = \phi_0(p)$ when $\epsilon = 0$ and then, by the expansion of the function $\Gamma_p$ and the estimations on $(\Lambda_\epsilon,w_\epsilon,\varphi_\epsilon)$, we have that when $\epsilon = 0$ the function $\hat \phi_{\epsilon,0}(y)$ is equal to $\phi_1(y)$.
In a neighborhood of $\partial \mathring B_1$ the metric $\hat g$ converge, for $\epsilon = 0$, to the euclidean metric, and from this it follows that $F(p,0,0)$ is the normal derivative of $\phi_1$ at $\partial \mathring B_1$ minus its euclidean mean, hence equal to 0.

\medskip

Our next step is to compute $F(p,0,s\bar w)$, and for this we have to study $F(p,\epsilon,s\, \bar w)$. Writing $\bar v = s \, \bar w$, we can consider a parameterization $Y$ of $B_{2\epsilon}^{g}(p)$ given by the following expression :
\[
Y (y) : = \mbox{Exp}_p^{\bar g} \left( \left(1 + \chi_1(y)\, v_0 + s\, \chi_2 ( y) \, \left(\bar w \left(\frac{y}{|y|} \right) \right) \right) \, \sum_i y^i \, E_i \right) 
\]
where $\bar g$ is the dilated metric $\epsilon^{-2}\, g$, $y$ belongs to the euclidean ball $\mathring B_2$ of radius 2 centered at 0, $\chi_1$ is a cutoff function identically equal to $1$ when $0 < |y| \leq 4/3$ and identically equal to $0$ when $5/3 \leq |y| \leq 2$, $\chi_2$ is a cutoff function identically equal to $1$ when $3/4 \leq |y| \leq 4/3$ and identically equal to $0$ when $0 <|y| \leq 1/2$ and $5/3 \leq |y| \leq 2$, and $v_0 = v_0(p, \epsilon, s\, \bar w)$. We set 
\[
\hat g : = Y^* \bar  g.
\]
over $\mathring B_2$. It is an extension of the metric $\hat g$ that we defined before on $\mathring B_1$.
We remark that $\hat \phi_{\epsilon,0}  : = Y^* \phi_{\epsilon,0}$ is a solution on $\mathring B_2$ of
\[
\Delta_{\hat{g}} \, \hat \phi_{\epsilon,0} + \hat \lambda_{\epsilon,0} \, \hat \phi_{\epsilon,0} = 0 
\]
where $\hat \lambda_{\epsilon,0} =\bar \lambda_{\epsilon,0} = \epsilon^{2}\, \lambda_{\epsilon,0}$. If we set $\bar \phi_{\epsilon,0}(y) = \phi_{\epsilon,0}(\epsilon\,y)$ in a neighborhood of $\partial \mathring B_1$, where $x=\epsilon\,y$ are the normal coordinates near $p$ that we defined in paragraph 3, we have
\begin{equation}\label{phialbordo}
\hat \phi_{\epsilon,0} (y) = \bar \phi_{\epsilon,0}  ((1+ v_0 + s \, \bar w(y) ) \, y ) \, ,
\end{equation}
on $\partial \mathring B_{1}$. Writing the first eigenfunction of $-\Delta_{\bar g}$ on $B^{\bar g}_{1+v}(p)$ as $ \phi = \phi_{\epsilon,0} + \psi$ and $\bar \lambda = \bar \lambda_{\epsilon,0} + \tau$, we find that 
\begin{equation}
\label{f001}
\left\{
\begin{array}{rcccl}
	(\Delta_{\bar{g}} + \bar \lambda_{\epsilon,0}) \, \psi + \tau \, \psi + \tau \, \phi_{\epsilon,0} & = & 0 & \textnormal{in} & M \setminus B_{1+v}^{\bar g}(p) \\[3mm]
 \psi & = & - \phi_{\epsilon,0} & \textnormal{on} & \partial B_{1+v}^{\bar g}(p) 
\end{array}
\right.
\end{equation}
where we can normalize as
\begin{equation}\label{integraleee}
\int_{M \setminus B^{\bar g}_{1+v}(p)} (\phi_{\epsilon,0} + \psi)^2 \textnormal{ dvol}_{\bar g} = \int_{M \setminus B^{\bar g}_{1+v_0}(p)} \phi_{\epsilon,0}^2 \textnormal{ dvol}_{\bar g} 
\end{equation}
(the $v_0$ in the second integral is evaluated at $\bar v =0$) and we have the condition on the volume of the domain
\begin{equation}
\label{formula-new-20}
\textnormal{ Vol}_{\hat g}(\mathring B_{1}) = \textnormal{ Vol}_{\mathring g} \, (\mathring B_1)
\end{equation}

Obviously $\psi$, $\tau$ and $v_0$ are smooth functions of $s$. When $s =0$, we have $ \phi =\phi_{\epsilon,0}$, $\bar \lambda = \bar \lambda_{\epsilon,0}$ and $v_0 =\mathcal{O}(\epsilon^2)$. Therefore,  $ \psi$ and $\tau$ vanish when $s=0$. We set 
\[
\dot \psi = \partial_s \psi |_{s=0} , \qquad \dot \tau = \partial_s \tau |_{s=0} , \qquad \mbox{and} \qquad \dot v_0 = \partial_s v_0 |_{s=0} ,
\]
Differentiating (\ref{f001}) with respect to $s$ and evaluating the result at $s=0$, we obtain 
\begin{equation}
\label{dotpsi}
\left\{
\begin{array}{rlllll}
	(\Delta_{\bar g} + \bar \lambda_{\epsilon,0})\, \dot \psi  + \dot \tau \, \phi_{\epsilon,0}  & = & 0 & \textnormal{in} & M \setminus B_{1+v_0}^{\bar g}(p) \\[3mm]
	\dot \psi & = & - \bar g(\nabla \phi_{\epsilon,0},\bar \nu) \, ( \dot v_0 + \bar w) & \textnormal{on} & \partial B_{1+v_0}^{\bar g}(p)
\end{array}
\right.
\end{equation}
where $v_0$ is evaluated at $s=0$. Observe that $\phi_{\epsilon,0}$ on $\partial B^{\bar g}_{1+v}(p)$ is equal to $\hat \phi_{\epsilon,0}$ on $\partial \mathring B_1$, then the second equation of (\ref{dotpsi}) follows from (\ref{phialbordo}). 

\medskip

Differentiating (\ref{integraleee}) with respect to $s$ and evaluating the result at $s=0$, we obtain that $\psi$ is $L^2$-orthogonal to $\phi_{\epsilon,0}$ on $B^{\bar g}_{1+v_0}(p)$. Hence 
\[
\phi = \phi_{\epsilon,0} + s \, \dot \psi + {\mathcal O} (s^2)
\]
where $\dot \psi$ is the solution of (\ref{dotpsi}) $L^2$-orthogonal to $\phi_{\epsilon,0}$. Differentiating (\ref{formula-new-20}) with respect to $s$ and evaluating the result at $s=0$, we obtain 
\begin{equation}
\label{formula-new-200}
\int_{S^{n-1}} ( \dot v_0 + \bar w) \, \textnormal{ dvol}_{\hat g}  =0
\end{equation}
where the metric $\hat g$ is evaluated at $s=0$.
Since the discrepancy between the metric $\hat g$ and the euclidean metric $\mathring g$ at $\partial \mathring B_1$ can be estimated by a constant times $\epsilon^2$ when $s=0$, and the euclidean average of $\bar w$ is $0$, we get that 
\[
\dot v_0 = {\mathcal O} (\epsilon^2)
\]
and then from the Taylor expansion of $v_0$ with respect to $s$ we get  
\[
v_0 =  {\mathcal O} (\epsilon^2) + {\mathcal O} (s^2)
\]
Now, in $\mathring B_{4/3} \setminus \mathring B_{1}$, we have 
\[
\begin{array}{rlll}
\hat \phi  (y) & = & \bar \phi_{\epsilon,0} \left(\big(1+ v_0(0) + s \, \bar w(y/|y|)\big) \, y \right) + s \, \dot \psi (y) + {\mathcal O} (s^2) \\[3mm]
                      & =  &  \bar\phi_{\epsilon,0} ((1+ v_0(0))\,y) + s\, \left( \hat g \left( \nabla \bar \phi_{\epsilon,0}((1+ v_0(0))\,y),(\dot v_0 +  \bar w (y/|y|))\,y\right)  + \dot \psi \right) + {\mathcal O} (s^2)
\end{array}
\]
where we denoted $v_0(p,\epsilon,0) = v_0|_{s=0} = v_0(0)$.
To complete the proof of the result, it suffices to compute the normal derivative of the function $\hat \phi$ when the normal is computed with respect to the metric $\hat g$. We use polar coordinates $y = r \, z$ where $r >0$ and $z \in S^{n-1}$. Then the metric $\hat g$ can be expanded in $\mathring B_{4/3}\setminus \mathring B_{3/4}$ as
\[
\hat g = (1 + v_0 + s \bar w )^2 \, dr^2  + 2 \, s \, (1+ v_0 + s \bar w) \, r \, d \bar w \, dr +  r^2 \,  (1+ v_0 + s \bar w)^2  \, \mathring h + s^2 \, r^2 \, d\bar w^2\, + {\mathcal O} (\epsilon^2)
\]
where $\mathring h$ is the metric on $S^{n-1}$ induced by the Euclidean metric. It follows from this expression, together with the estimation of $v_0$, that the unit normal vector field to $\partial \mathring B_1$ for the metric $\hat g$ is given by
\[
\hat \nu  = \left( ( 1 + s \, \bar w )^{-1} +  {\mathcal O} (s^2)\right)  \,  \partial_r + {\mathcal O} (s) \, \partial_{z_j}\, + {\mathcal O} (\epsilon^2)
\]
where $\partial_{z_j}$ are vector fields induced by a parameterization of $S^{n-1}$. Using this, we conclude that 
\begin{equation}\label{operepsilon}
\hat g ( \nabla  \hat \phi , \hat \nu ) = \partial_r \bar \phi_{\epsilon,0}(y) + {\mathcal O} (s) \, \partial_{z_j}\, \bar \phi_{\epsilon,0}(y) + s \, \left( \bar w \,  \partial_r^2\, \bar \phi_{\epsilon,0}(y) + \partial_r \dot \psi \right) + {\mathcal O}(s^2) \, + {\mathcal O} (\epsilon^2)
\end{equation}
on $\partial \mathring B_1$. When $\epsilon = 0$ we have that $\bar \phi_{\epsilon,0}(y) = \phi_1(y)$. It follows that $F(p,0,s\bar w)$, up to terms of the order ${\mathcal O}(s^2)$, is given by the term 
\[
\left. \partial_r \phi_1\right|_{\partial \mathring B_1} + s\,\bar w \left. \partial^2_r \phi_1\right|_{\partial \mathring B_1} + s\, \lim_{\epsilon \rightarrow 0} \left.\partial_r \dot \psi\right|_{\partial \mathring B_1}
\]
minus its euclidean mean, where the limit is understood in the pointwize sens. 
We need the following result.

\begin{Lemma}\label{penultimo}
Evaluate $v_0$ at $s=0$. Let $\nu \in (2-n, 0)$ if $n \geq 3$ and $\nu \in (0,1)$ if $n = 2$. Let $H_\varphi$ be the function defined in the paragraph 4. 
For all $\epsilon$ small enough there exist a constant $\dot \tau$ and  $(K_\epsilon,\varphi_\epsilon,\eta_\epsilon)$ in a neighborhood of $(0,0,0)$ in $\mathbb{R} \times \mathcal{C}^{2,\alpha}_m(S^{n-1}) \times \mathcal{C}^{2,\alpha}_\nu(M \setminus B_{1+v_0}^{\bar g}(p))$ such that the function
\begin{equation}
\dot \psi = K_\epsilon + \chi\, (\psi + H_{\varphi_\epsilon}) + \eta_\epsilon
\end{equation}
defined in $M \setminus B_{1+v_0}^{\bar g}(p)$, is the solution of (\ref{dotpsi}) $L^2$-orthogonal to $\phi_{\epsilon,0}$, where $\chi$ is a cut-off function equal to 1 in $B^{\bar g}_{R_0/\epsilon}(p)$ and equal to 0 out of $B^{\bar g}_{2R_0/\epsilon}(p)$ and $\psi$ is defined by (\ref{c00}).
Moreover the following estimations hold :
	\[
|K_\epsilon| \leq c\, \left(\epsilon^2 + \epsilon^{n-1}\right) \qquad \textnormal{and} \qquad \left\|\varphi_\epsilon\right\|_{L^{\infty}(S^{n-1})} \leq c\, \left(\epsilon^2 + \epsilon^{n-1}\right) 
\]
\[
\textnormal{and} \qquad \| \eta_\epsilon \|_{\mathcal{C}^{2,\alpha}_\nu(M \setminus B_{1+v_0}^{\bar g}(p))}\leq c\, \left(\epsilon^2 + \epsilon^{n-1}\right) 
\]
	\end{Lemma}

{\bf Proof.} 
Let us choose $\dot \psi$ in the form
\begin{equation}\label{dotpsichosen}
\dot \psi =K + \chi\, (\psi + H_{\varphi}) + \eta
\end{equation}
for some $(K, \varphi,\eta) \in \mathbb{R} \times \mathcal{C}^{2,\alpha}_m(S^{n-1}) \times \mathcal{C}^{2,\alpha}(M \setminus B_{1+v_0}^{\bar g}(p))$, where $\chi$ is a cut-off function equal to 1 in $B^{\bar g}_{R_0/\epsilon}(p)$ and equal to 0 out of $B^{\bar g}_{2R_0/\epsilon}(p)$.
Then $\dot \psi$ satisfy the first equation of (\ref{dotpsi}), if and only if :
\begin{equation}\label{eta}
\begin{array}{lll}
\left( \Delta_{\bar  g} + \bar  \lambda_{\epsilon,0}\right)\, \eta &= & - \psi \, \Delta_{\bar  g} \chi - \chi\, \Delta_{\bar  g} \psi -2\, \nabla^{\bar  g} \psi\, \nabla^{\bar  g} \chi - H_{\varphi} \, \Delta_{\bar  g} \chi - \chi\, \Delta_{\bar  g} H_{\varphi} -2\, \nabla^{\bar  g} H_{\varphi}\, \nabla^{\bar  g} \chi \\[3mm]
& & \qquad - \bar  \lambda_{\epsilon,0}\, \chi\, (\psi + H_{\varphi}) - \bar  \lambda_{\epsilon,0}\, K -\dot \tau \, \phi_{\epsilon,0} 
\end{array}
\end{equation}
For $n\geq 3$ and $\nu \in (2-n,0)$, the operator
\[
\big(\Delta_{\bar g} + \bar \lambda_{\epsilon,0} \big): \mathcal{C}^{2,\alpha}_{\nu,\bot,0}(M \setminus B_{1+v_0}^{\bar g}(p)) \longrightarrow \mathcal{C}^{0,\alpha}_{\nu-2,\bot}(M \setminus B_{1+v_0}^{\bar g}(p)),
\]
where the subscript $\bot$ is meant to point out that functions are $L^2$-orthogonal to $\phi_{\epsilon,0}$, and the subscript $0$ is meant to point out that functions satisfy the 0 Dirichlet (CASE 1) or 0 Neumann (CASE 2) condition on $\partial M$ and the 0 Dirichlet condition on $\partial B_{1+v_0}^{\bar g}(p)$, is an isomorphism. For $n=2$ and $\nu \in (0,1)$ the same result holds for the operator 
\[
\big(\Delta_{\bar g} + \bar \lambda_{\epsilon,0} \big): \tilde \chi\, \mathbb{R} \oplus \mathcal{C}^{2,\alpha}_{\nu,\bot,0}(M \setminus B_{1+v_0}^{\bar g}(p)) \longrightarrow \mathcal{C}^{0,\alpha}_{\nu-2,\bot}(M \setminus B_{1+v_0}^{\bar g}(p)).
\]
where $\tilde \chi$ is a cutoff function equal to 1 in a neighborhood of the origin. The proof of this facts has been given in paragraph 4.

\medskip 

To semplify the notation let us define
\[
\begin{array}{lll}
A &:=& - \psi \, \Delta_{\bar  g} \chi -2\, \nabla^{\bar  g} \psi\, \nabla^{\bar  g} \chi - H_{\varphi} \, \Delta_{\bar  g} \chi -2\, \nabla^{\bar  g} H_{\varphi}\, \nabla^{\bar  g} \chi \\[3mm]
B &:=& - \chi\, \Delta_{\bar  g} \psi - \bar  \lambda_{\epsilon,0}\, \chi\, (\psi + H_{\varphi})  - \chi\, \Delta_{\bar  g} H_{\varphi} - \bar  \lambda_{\epsilon,0}\, K\\[3mm]
C &:=& - \dot \tau\,  \phi_{\epsilon,0} 
\end{array}
\]
Equation (\ref{eta}) becomes
\[
( \Delta_{ \bar  g}+ \bar  \lambda_{\epsilon,0})\, \eta = A + B +C
\]
By the last result, if we chose $\dot \tau$ in order to verify
\begin{equation}\label{ABCDbis}
\int_{M \setminus B_{1+v_0}^{\bar g}(p)} (A + B +C)\, \phi_{\epsilon,0}
 = 0
\end{equation}
there exists a solution $\eta = \eta(\epsilon,K,\varphi) \in \mathcal{C}^{2,\alpha}_{\nu,\bot,0}(M \setminus B_{1+v_0}^{\bar g}(p))$ (or $\tilde \chi\, \mathbb{R} \oplus \mathcal{C}^{2,\alpha}_{\nu,\bot,0}(M \setminus B_{1+v_0}^{\bar g}(p))$ if $n=2$) to equation (\ref{eta}) for all $\epsilon$ small enough, for all constant $K$ and all function $\varphi$, and then
\[
\dot \psi = K + \chi\, (\psi + H_{\varphi}) + \eta
\]
satisfy the first equation of (\ref{dotpsi}). 

\medskip

We want now to give some estimations on the function $\eta$. By the previous results and Lemma \ref{lemmaharmonic} we have the following estimations :
\medskip

\begin{itemize}
	\item $\| A \|_{\mathcal{C}^{0,\alpha}_{\nu-2}(M \setminus B_{1+v_0}^{\bar g}(p))} \leq c\, \epsilon^{n-1}\, \left(1 + \|\varphi\|_{L^{\infty}(S^{n-1})}\right)$

	\medskip
	
	\item $\| B \|_{\mathcal{C}^{0,\alpha}_{\nu-2}(M \setminus B_{1+v_0}^{\bar g}(p)} \leq c\, \epsilon^{2}\, \left( 1 + \|\varphi\|_{L^{\infty}(S^{n-1})}\right)$
	\medskip
\end{itemize}
In particular we get that 
\[
\dot \tau \leq c\, \left(\epsilon^2 + \epsilon^{n-1}\right)\, \left(1 + \|\varphi\|_{L^{\infty}(S^{n-1})}\right)
\]
and then
\[
\| A + B+ C \|_{\mathcal{C}^{0,\alpha}_{\nu-2}(M \setminus B_{1+v_0}^{\bar g}(p)))}  \leq c\, \left(\epsilon^2 + \epsilon^{n-1}\right)\, \left(1 + \|\varphi\|_{L^{\infty}(S^{n-1})}\right)
\]
This give us an estimation on the function $\eta$ that we found before:
\[
\| \eta \|_{\mathcal{C}^{2,\alpha}_{\nu}(M \setminus B_{1+v_0}^{\bar g}(p))}  \leq c\, \left(\epsilon^2 + \epsilon^{n-1}\right)\, \left(1 + \|\varphi\|_{L^{\infty}(S^{n-1})}\right).
\]
Summarizing, we have proved the following : For all $\varphi \in \mathcal{C}^{2,\alpha}_{m}(S^{n-1})$, for all constant $K$, for all $\epsilon$ small enough, there exists a function $\eta(\epsilon,K,\varphi) \in \mathcal{C}^{2,\alpha}_{\nu,\bot, 0}(M \setminus B_{1+v_0}^{\bar g}(p))$ such that (\ref{dotpsichosen}) is a positive solution of the first equation of (\ref{dotpsi}). Moreover there exists a positive constant $c$ such that
\[
\| \eta \|_{\mathcal{C}^{2,\alpha}_{\nu}(M \setminus B_{1+v_0}^{\bar g}(p))}  \leq c\, \left(\epsilon^2 + \epsilon^{n-1}\right)\, \left(1 + \|\varphi\|_{L^{\infty}(S^{n-1})}\right).
\]
\medskip

Now we have to make attention to the second equation of (\ref{dotpsi}). Let us define
\[
Z(\epsilon,K,\varphi) := \big [K + \chi(y)\, (\psi(y) + H_{\varphi}(y)) + \eta(\epsilon,\varphi)(y)\big]_{y \in S^{n-1}}.
\]
We remark that $Z$, that represents the boundary value of the solution of the first of (\ref{dotpsi}) we found above, is well defined in a heighborhood of $(0,0,0)$ in $(0,+\infty) \times \mathbb{R} \times \mathcal{C}^{2,\alpha}_m(S^{n-1})$, and takes its values in $\mathcal{C}^{2,\alpha}(S^{n-1})$. It is easy to compute the differential of $Z$ with respect to $K$ and $\varphi$ at $(0,0,0)$ :
\[
\big(\partial_{\varphi} Z(0,0,0)\big)(\tilde K) = \tilde K.
\]
\[
\big(\partial_{\varphi} Z(0,0,0)\big)(\tilde \varphi) = \tilde \varphi.
\]
We can estimate $Z(\epsilon,0,0)$ :
\[
\left\|Z(\epsilon,0, 0) + \partial_r \phi_1\, \bar w\right\|_{L^{\infty}(S^{n-1})} \leq c\, \left(\epsilon^2 + \epsilon^{n-1}\right) 
\]
and then
\[
\left\|Z(\epsilon,0, 0) + \bar g (\nabla \phi_{\epsilon,0}, \bar \nu)\, (\dot v_0 + \bar w) \right\|_{L^{\infty}(S^{n-1})} \leq c\, \left(\epsilon^2 + \epsilon^{n-1}\right) 
\]
The implicit function theorem applies to give the following : Let $\epsilon$ be small enough; then there exist $(K_\epsilon,\varphi_\epsilon)$ in a neighborhood of $(0,0)$ in $\mathbb{R} \times \mathcal{C}^{2,\alpha}_m(S^{n-1})$ such that 
(\ref{dotpsichosen}) is a positive solution of (\ref{dotpsi}). Moreover the following estimations hold :
	\[
|K_\epsilon| \leq c\, \left(\epsilon^2 + \epsilon^{n-1}\right) \qquad \textnormal{and} \qquad \left\|\varphi_\epsilon\right\|_{L^{\infty}(S^{n-1})} \leq c \,\left(\epsilon^2 + \epsilon^{n-1}\right)
	\]
	
\medskip 

Summarizing, we get the following existence result: for all $\epsilon$ small enough there exist a constant $\dot \tau$ and $(K_\epsilon, \varphi_\epsilon,\eta_\epsilon)$ in a neighborhood of $(0,0,0)$ in $\mathbb{R} \times \mathcal{C}^{2,\alpha}_m(S^{n-1}) \times \mathcal{C}^{2,\alpha}_\nu(M \setminus B_{1+v_0}^{\bar g}(p))$ (or $\mathbb{R} \times \mathcal{C}^{2,\alpha}_m(S^{n-1}) \times \tilde \chi\, \mathbb{R} \oplus \mathcal{C}^{2,\alpha}_\nu(M \setminus \{p\})$ if $n=2$) such that the function
\[
\dot \psi = K_\epsilon + \chi\, (\psi + H_{\varphi_\epsilon}) + \eta_\epsilon
\]
defined in $M \setminus B_{1+v_0}^{\bar g}(p)$, is solution of (\ref{dotpsi})
Moreover :
	\[
|K_\epsilon| \leq c\, \left(\epsilon^2 + \epsilon^{n-1}\right) \qquad \textnormal{and} \qquad \left\|\varphi_\epsilon\right\|_{L^{\infty}(S^{n-1})} \leq c\, \left(\epsilon^2 + \epsilon^{n-1}\right) 
\]
\[
\textnormal{and} \qquad \| \eta_\epsilon \|_{\mathcal{C}^{2,\alpha}_\nu(M \setminus B_{1+v_0}^{\bar g}(p))}\leq c\, \left(\epsilon^2 + \epsilon^{n-1}\right) 
\]
The last norm is over $\tilde \chi\, \mathbb{R} \oplus \mathcal{C}^{2,\alpha}_\nu(M \setminus \{p\})$ if $n=2$. This completes the proof of the result. \hfill $\Box$

\medskip

Using the previous lemma, we have that for $\epsilon$ small enough
\[
\left.\partial_r \dot \psi\right|_{\partial \mathring B_1} = \left.\partial_r \psi\right|_{\partial \mathring B_1} + \mathcal{O}(\epsilon^2)
\]
for $n \geq 3$ and
\[
\left.\partial_r \dot \psi\right|_{\partial \mathring B_1} = \left.\partial_r \psi\right|_{\partial \mathring B_1} + \mathcal{O}(\epsilon)
\]
for $n = 2$, because the estimation of $\eta_\epsilon$ is given on the weighted Holder spaces.
The statement of the Proposition \ref{l0} then follows at once from the fact that $\partial_r \phi_1$ is constant  while the term $ \bar w  \,  \partial_r^2 \phi_1  + \partial_r \psi $ has mean $0$ on the boundary $\partial \mathring B_1$. This completes the proof of the proposition.\hfill  $\Box$

\medskip

Now we denote by $L_\epsilon$ the linearization of $F$ with respect to $\bar v$, computed at the point $(p, \epsilon, 0)$. It is easy to check the~:
\begin{Lemma}There exists a constant $c >0$ such that, for all $\epsilon >0$ small enough we have the estimate
\[
\begin{array}{ll}
\| (L_\epsilon  - L_0) \, \bar v   \|_{\mathcal C^{1, \alpha}}\leq c \, \epsilon \, \| \bar v   \|_{\mathcal C^{2, \alpha}} & \textnormal{in the CASE 1 and}\,\,n \geq 3\\[3mm]
\| (L_\epsilon  - L_0) \, \bar v   \|_{\mathcal C^{1, \alpha}}\leq c \, \epsilon\, \log \epsilon \, \| \bar v   \|_{\mathcal C^{2, \alpha}} & \textnormal{in the CASE 1 and}\,\,n = 2\\[3mm]
\| (L_\epsilon  - L_0) \, \bar v   \|_{\mathcal C^{1, \alpha}}\leq c \, \epsilon^2\, \| \bar v   \|_{\mathcal C^{2, \alpha}} & \textnormal{in the CASE 2 and}\,\,n \geq 5\\[3mm]
\| (L_\epsilon  - L_0) \, \bar v   \|_{\mathcal C^{1, \alpha}}\leq c \, \epsilon^2\, \log \epsilon \, \| \bar v   \|_{\mathcal C^{2, \alpha}} & \textnormal{in the CASE 2 and}\,\,n = 4
\end{array}
\]
\label{le:4444}
\end{Lemma}
{\bf Proof :}  Clearly both $L_\epsilon$ and $L_0$ are first order differential operators. We already know the expression of $L_0$. We have
\[
L_\epsilon (\bar w) = \lim_{s \rightarrow 0} \frac{F(p,\epsilon, s\, \bar w) - F(p,\epsilon,0)}{s} .
\]
$F(p,\epsilon, s\, \bar w)$ is given by (\ref{operepsilon}) minus its mean, in the metric $\hat g$. $F(p,\epsilon,0)$, up to terms of order $\mathcal{O}(\epsilon^2)$, is given by $\partial_r \bar \phi_{\epsilon,0}(y)$ at $\partial \mathring B_1$ minus its the mean, in the metric $\hat g$ evaluated at $s=0$. The proof of the Lemma follows at once from Proposition \ref{le:3.3} and Lemma \ref{penultimo}.
\hfill $\Box$

\section{The proof of the main result}

We shall now prove that, for $\epsilon>0$ small enough, it is possible to solve the equation 
\[
F  (p, \epsilon, \bar v ) =0
\]
Unfortunately, we will not be able to solve this equation at once. Instead, we first prove the~:
\begin{Proposition}
There exists $\epsilon_0 >0$ such that, for all $\epsilon \in [0, \epsilon_0]$ and for all $p \in M$, there exists a unique function $\bar v = \bar v(p, \epsilon)$ and a vector $a  = a (p, \epsilon) \in \mathbb R^n$ such that 
\[
F (p, \epsilon, \bar v )  + \mathring g ( a , \cdot) = 0 
\]
The function $\bar v$  and the vector $a$ depend smoothly on $p$ and $\epsilon$ and we have 
\[
\begin{array}{ll}
|a | + \| \bar v \|_{\mathcal C^{2, \alpha}(S^{n-1})} \leq c \, \epsilon & \textnormal{in the CASE 1 and}\,\,n \geq 3\\[3mm]
|a | + \| \bar v \|_{\mathcal C^{2, \alpha}(S^{n-1})} \leq c \, \epsilon\, \log \epsilon & \textnormal{in the CASE 1 and}\,\,n = 2\\[3mm]
|a | + \| \bar v \|_{\mathcal C^{2, \alpha}(S^{n-1})} \leq c \, \epsilon^2 & \textnormal{in the CASE 2 and}\,\,n \geq 5\\[3mm]
|a | + \| \bar v \|_{\mathcal C^{2, \alpha}(S^{n-1})} \leq c \, \epsilon^2\, \log \epsilon & \textnormal{in the CASE 2 and}\,\,n =4
\end{array}
\]
\end{Proposition}
{\bf Proof :} We fix $p \in M$ and define
\[
\bar F (p, \epsilon, \bar v , a) : = F (p, \epsilon, \bar v )  + \mathring g (a, \cdot)
\]
It is easy to check that $\bar F$ is a smooth map from a neighborhood of ${(p,0,0,0)}$ in $M \times [0, \infty) \times \mathcal C^{2, \alpha}_m (S^{n-1}) \times \mathbb R^n$ into a neighborhood of $0$ in $\mathcal C^{1, \alpha} (S^{n-1})$. Moreover, 
\[
\bar F (p, 0, 0 , 0) =0
\]
and the differential of $\bar F$ with respect to $\bar v$, computed at ${(p,0,0,0)}$ is given by $H$. Finally the image of the linear map $a \longmapsto \mathring g (a , \cdot)$ is just the vector space $V_1$.  Thanks to the result of Proposition~\ref{H}, the implicit function theorem applies to get the existence of $\bar v$ and $a$, smoothly depending on $p$ and $\epsilon$ such that $\bar F (p, \epsilon, \bar v,a  ) = 0$. The estimates for $\bar v$ and $a$ follow at once from Proposition~\ref{le:3.3}. This completes the proof of the result. \hfill $\Box$

\medskip

In view of the result of the previous Proposition, it is enough to show that, provided that $\epsilon$ is small enough, it is possible to choose the point $p \in M$ such that $a (p,\epsilon) =0$.  We claim that, there exists a constant $\bar C >0$ (only depending on $n$) such that
\[
\begin{array}{ll}
\Theta (a(p, \epsilon ))  =  -\epsilon\,  \tilde C \, \nabla^g \phi_0 (p) + \mathcal O (\epsilon^2) & \textnormal{in the CASE 1 and}\,\,n \geq 3\\[3mm]
\Theta (a(p, \epsilon ))  =  - \epsilon\, \log \epsilon \, \tilde C \, \nabla^g \phi_0 (p) + \mathcal O (\epsilon)& \textnormal{in the CASE 1 and}\,\,n = 2\\[3mm]
\Theta (a(p, \epsilon ))  = - \epsilon^3  \, \tilde C \, \nabla^g \textnormal{ Scal} (p) + \mathcal O (\epsilon^4) & \textnormal{in the CASE 2 and}\, \, n \geq 5\\[3mm]
\Theta (a(p, \epsilon ))  = - \epsilon^3\, \log\epsilon  \, \tilde C \, \nabla^g \textnormal{ Scal} (p) + \mathcal O (\epsilon^3) & \textnormal{in the CASE 2 and}\,\, n= 4 
\end{array}
\]

\medskip

For  all $b \in \mathbb R^n$ we compute 
\[
\begin{array}{rlllll}
\displaystyle \int_{S^{n-1}} \, \mathring g(a , \cdot ) \, \mathring g(b, \cdot ) \, \textnormal{ dvol}_{\mathring g} & = &\displaystyle  - \int_{S^{n-1}} \, F (p, \epsilon, \bar v) \, \mathring g(b, \cdot ) \, \textnormal{ dvol}_{\mathring g}\\[3mm]
& = &  - \displaystyle \int_{S^{n-1}} \, ( F (p, \epsilon, 0)  + L_0 \bar v ) \, \mathring g(b, \cdot ) \, \textnormal{ dvol}_{\mathring g}\\[3mm]
& & -  \displaystyle \int_{S^{n-1}} \, (F (p, \epsilon, \bar v)  - F(p, \epsilon, 0)  - L_\epsilon \bar v) \, \mathring g(b, \cdot ) \, \textnormal{ dvol}_{\mathring g}\\[3mm]
& &  -  \displaystyle \int_{S^{n-1}} \, (L_\epsilon  - L_0) \bar v \, \mathring g(b, \cdot ) \, \textnormal{ dvol}_{\mathring g}\\[3mm]
\end{array}
\]

Now, we use the fact that $\bar v$ is $L^2(S^{n-1})$-orthogonal to linear functions and hence so is $L_0 \, \bar v$. Therefore, 
\[
{\displaystyle \int_{S^{n-1}} \,  L_0 \, \bar v  \, \mathring g(b, \cdot ) \, \textnormal{ dvol}_{\mathring g} = 0}
\] 
Using the fact that 
\[
\begin{array}{ll}
\bar v = \mathcal O (\epsilon) & \textnormal{in the CASE 1 and}\,\,n \geq 3\\[3mm]
\bar v = \mathcal O (\epsilon\, \log\epsilon) & \textnormal{in the CASE 1 and}\,\,n = 2\\[3mm]
\bar v = \mathcal O (\epsilon^2) & \textnormal{in the CASE 2 and}\,\,n \geq 5\\[3mm]
\bar v = \mathcal O (\epsilon^2\, \log\epsilon) & \textnormal{in the CASE 2 and}\,\,n = 4
\end{array}
\]
we get 
\[
\begin{array}{ll}
F (p, \epsilon, \bar v)  - F(p, \epsilon, 0)  - L_\epsilon \bar v =  \mathcal O (\epsilon^2)  & \textnormal{in the CASE 1 and}\,\,n \geq 3\\[3mm]
F (p, \epsilon, \bar v)  - F(p, \epsilon, 0)  - L_\epsilon \bar v =  \mathcal O (\epsilon^2\, (\log\epsilon)^2)  & \textnormal{in the CASE 1 and}\,\,n = 2\\[3mm]
F (p, \epsilon, \bar v)  - F(p, \epsilon, 0)  - L_\epsilon \bar v =  \mathcal O (\epsilon^4)  & \textnormal{in the CASE 2 and}\,\,n \geq 5\\[3mm]
F (p, \epsilon, \bar v)  - F(p, \epsilon, 0)  - L_\epsilon \bar v =  \mathcal O (\epsilon^4\, (\log\epsilon)^2)  & \textnormal{in the CASE 2 and}\,\,n = 4 
\end{array}
\]
Similarly, it follows from the result of Proposition~\ref{le:4444} that
\[
\begin{array}{ll}
( L_\epsilon  - L_0) \,  \bar v = \mathcal O (\epsilon^2) & \textnormal{in the CASE 1 and}\,\,n \geq 3\\[3mm]
( L_\epsilon  - L_0) \,  \bar v = \mathcal O (\epsilon^2\, (\log\epsilon)^2) & \textnormal{in the CASE 1 and}\,\,n = 2\\[3mm]
( L_\epsilon  - L_0) \,  \bar v = \mathcal O (\epsilon^4) & \textnormal{in the CASE 2 and}\,\,n \geq 5\\[3mm]
( L_\epsilon  - L_0) \,  \bar v = \mathcal O (\epsilon^4\, (\log\epsilon)^2) & \textnormal{in the CASE 2 and}\,\,n = 4 
\end{array}
\]
The claim then follows from the estimates in Proposition~\ref{le:3.3} {and the fact that
\[
\int_{S^{n-1}} \, \mathring g(a , \cdot ) \, \mathring g(b, \cdot ) \, \textnormal{ dvol}_{\mathring g} = g\big(\Theta(a),\Theta(b)\big) \,  \int_{S^{n-1}} (x_1)^2 \, \textnormal{ dvol}_{\mathring g} = \frac{1}{n}\, \textnormal{ Vol}_{\mathring g} (S^{n-1})\, g\big(\Theta(a),\Theta(b)\big).
\]}

\medskip
 
Now if we assume that  $p_0$ is a nondegenerate critical point of the function $\phi_0$ (CASE 1) or a nondegenerate critical point of the scalar curvature function (CASE 2), we can apply once more the implicit function theorem to solve the equations 
\[
\begin{array}{ll}
G(\epsilon, p) : = \epsilon^{-1}\, \Theta (a(p, \epsilon) ) = 0 & \textnormal{in the CASE 1 and}\,\,n \geq 3\\[3mm]
G(\epsilon, p) : = (\epsilon\, \log\epsilon)^{-1}\, \Theta (a(p, \epsilon) ) = 0 & \textnormal{in the CASE 1 and}\,\,n = 2\\[3mm]
G(\epsilon, p) : = \epsilon^{-3}\, \Theta (a(p, \epsilon) ) = 0 & \textnormal{in the CASE 2 and}\, \, n \geq 5\\[3mm]
G(\epsilon, p) : = \epsilon^{-3}\, (\log\epsilon)^{-1}\, \Theta (a(p, \epsilon) ) = 0 & \textnormal{in the CASE 2 and}\,\, n= 4 
\end{array}
\]
It should be clear that $G$ depends smoothly on $\epsilon \in [0, \epsilon_0)$ and $p \in M$. Moreover we have
\[
G(0, p) =  - \tilde C \, \nabla^g \phi_0(p) 
\]
in the CASE 1 and 
\[
G(0, p) =  - \tilde C \, \nabla^g \textnormal{ Scal} (p) 
\]
in the CASE 2. Hence, under the hypothesis on $p_0$ of the main theorem, $G(0, p_0) =0$ in both cases. By assumption the differential of $G$ with respect to $p$, computed at $p_0$ is invertible. Therefore, for all $\epsilon$ small enough there exists $p_\epsilon$ close to $p_0$ such that 
\[
\Theta (a(p_\epsilon, \epsilon) ) = 0
\]
In addition we have 
\[
\mbox{dist} (p_0 , p_\epsilon ) \leq c \, \epsilon
\]
This completes the proof the Theorem~\ref{maintheorem}. 

\section{Appendix}

We recall here some results demonstrated in [\ref{Pac-Sic}].
\begin{Lemma}
For all $\sigma =1, \ldots, n$, we have
$$
\sum_{i,j, k, \ell,m} \int_{S^{n-1}} R_{ikj\ell,m} \, x^{i} \, x^{j} \, x^{k} \, x^{\ell} \, x^{m} \, x^{\sigma} \, \textnormal{ dvol}_{\mathring g} = 0.
$$
\end{Lemma}

\medskip

\begin{Lemma}
For all $\sigma =1, \ldots, n$, we have
$$
\sum_{j,k, \ell} \int_{S^{n-1}} R_{\cdot kj\ell ,\cdot} \, x^{j} \, x^k \, x^{\ell} \,  x^{\sigma} \, \textnormal{ dvol}_{\mathring g} = 0.
$$
\end{Lemma}

\medskip

\begin{Lemma}
For all $\sigma =1, \ldots, n$, we have
$$
\sum_{i, \ell ,m} \, \int_{S^{n-1}} \, R_{i\ell ,m}  \, x^{i} \, x^{\ell} \, x^m \, x^{\sigma} \, \textnormal{ dvol}_{\mathring g} = \frac{2}{n(n+2)} \, \textnormal{ Vol}_{\mathring g} (S^{n-1})  \,\textnormal{ Scal}_{,\sigma} 
$$
\end{Lemma}

\medskip

Moreover we prove the

\begin{Lemma}
For all $\sigma =1, \ldots, n$, we have
$$
\sum_{t} \, \int_{S^{n-1}} \, {\textnormal{ Scal}}_{,t}\, x^{t} \, x^{\sigma} \, {\textnormal{ dvol}}_{\mathring g} = \frac{1}{n} \, \textnormal{ Vol}_{\mathring g} (S^{n-1})  \,{\textnormal{ Scal}}_{,\sigma} 
$$
\end{Lemma}

{\bf Proof :} We find that $\displaystyle \int_{S^{n-1}} \, {\textnormal{Scal}}_{,t}\, x^{t} \, x^{\sigma} \, \textnormal{ dvol}_{\mathring g} = 0$ unless the indices $t$ and $\sigma$ are equal. Then
\[
\sum_{t} \, \int_{S^{n-1}} \, {\textnormal{Scal}}_{,t}\, x^{t} \, x^{\sigma} \, {\textnormal{ dvol}}_{\mathring g} = {\textnormal{ Scal}}_{,\sigma} \int_{S^{n-1}} \, (x^{\sigma})^2 \, {\textnormal{ dvol}}_{\mathring g} = \frac{1}{n} \, {\textnormal{ Vol}}_{\mathring g} (S^{n-1})  \,\textnormal{ Scal}_{,\sigma} 
\]

\end{document}